\documentclass[a4paper,10pt]{article}
\usepackage{pkg}
\usepackage{ncmd}

\title{Rank of Matrices Arising out of Singular Kernel Functions\thanks{Submitted to the editors: \today}}
\author[4]{Sumit Singh\,\orcidlink{0009-0002-5581-5349}\thanks{\url{sumit1315singh@gmail.com}, \url{ma22d027@smail.iitm.ac.in}}}
\author[1,2,3,4]{Sivaram Ambikasaran\,\orcidlink{0000-0003-2978-6281}\thanks{\url{sivaambi@dsai.iitm.ac.in}, \url{sivaambi@alumni.stanford.edu}}}

\affil[1]{Wadhwani School of Data Science and Artificial Intelligence}
\affil[2]{Robert Bosch Centre for Data Science and Artificial Intelligence}
\affil[3]{Department of Data Science and Artificial Intelligence, IIT Madras, Chennai, India}
\affil[4]{Department of Mathematics, IIT Madras, Chennai, India}

\date{}

\begin{document}

\maketitle

\begin{abstract}
    Kernel functions are frequently encountered in differential equations and machine learning applications. In this work, we study the rank of matrices arising out of the kernel function $K: X \times Y \mapsto \mathbb{R}$, where the sets $X, Y \in \mathbb{R}^d$ are hypercubes that share a boundary. The main contribution of this work is the analysis of rank of such matrices where the particles (sources/targets) are \textbf{\emph{arbitrarily distributed}} within these hypercubes.
    To our knowledge, this is the first work to formally investigate the rank of such matrices for \textbf{\emph{arbitrary distribution of particles}}. We model the arbitrary distribution of particles to arise from an underlying random distribution and obtain bounds on the expected rank and variance of the rank of the kernel matrix corresponding to various neighbor interactions. These bounds are useful for understanding the performance and complexity of hierarchical matrix algorithms (especially hierarchical matrices satisfying the weak-admissibility criterion) for an arbitrary distribution of particles.
     We also present numerical experiments in one-, two-, and three-dimensions, showing the expected rank growth and variance of the rank for different types of interactions. The numerical results, not surprisingly, align with our theoretical predictions.
\end{abstract}

\noindent\textbf{Keywords:}  Probability distributions, Numerical rank, $n$-body problems, Hierarchical matrices,  Low-rank matrix approximation, Central Limit Theorem, Normal approximation. \\

\noindent\textbf{AMS Subject Classifications:} 65F55, 65D40, 65D12. 

\section{Introduction} \label{sec: Introduction}

 
    
Recent years have seen significant strides in the development of \textbf{matrices arising out of kernel functions} frequently occur in areas such as partial differential equations (PDEs) \cite{greengard1987fast} \cite{massei2022hierarchical} \cite{dick2013high} \cite{fornberg2015solving}, integral equations \cite{ho2016hierarchical}, inverse problems \cite{ambikasaran2013latge}\cite{harlim2020kernel}, Gaussian processes \cite{ambikasaran2015fast} \cite{nowak2013kriging}, graph theory \cite{kriege2020survey} \cite{vishwanathan2010graph} and kernel methods for addressing many complex machine learning and data analysis tasks \cite{cortes1995support} \cite{gray2000n} \cite{huang2006extreme}. Despite their wide applicability, these matrices \textit{are often large and dense}, as the underlying kernel functions are not compactly supported. This makes storage and computation of matrix operations (such as matrix-vector products, solving linear systems, and matrix-factorization, etc.) computationally expensive and memory intensive.

Despite these challenges, kernel matrices exhibit low-rank structure that can be exploited to overcome these issues. We can significantly reduce the storage requirements and accelerate computational processes by leveraging their low-rank approximations. The literature on exploiting rank-structuredness is extensive, and we refer interested readers to the works \cite{drineas2005nystrom} \cite{yokota2017fast} \cite{borm2003introduction} \cite{sarlos2006improved} \cite{halko2011finding}  and the references therein for an in-depth review.  This approach not only reduces computational costs but also allows for more accurate solutions to complex problems in scientific computing, engineering, and data science. More specifically, exploiting rank-structuredness is very useful in fields such as machine learning and signal processing, where it helps in optimized data compression, noise reduction, and the extraction of meaningful insights from large datasets.

One of the most frequently encountered rank-structured matrices arising out of $n$-body problems are \textbf{hierarchical low-rank} matrices. The initial work by Barnes and Hut \cite{barnes1986hierarchical}, which reduced the computational complexity to perform matrix-vector product from $\mclo{n^2}$ to $\mclo{n\log n}$. Greengard and Rokhlin in Fast Multipole Method (FMM) \cite{greengard1987fast} further reduced to $\mclo{n}$. FMM and Barnes-Hut algorithms leverage separable extension of kernel functions for far-away interactions. In terms of matrices, this corresponds to approximating sub-matrices corresponding to far-away interactions by a low-rank matrix. This interpretation generalises to hierarchical matrices. Depending on what sub-matrices are approximated by a low-rank matrix, we have different hierarchical structures. Some of the widely used such representations are Hierarchically Off-Diagonal Low-Rank (HODLR), Hierarchically Semi-Separable (HSS), and $\mathcal{H}$-matrix. For a detailed literature review on hierarchical matrices and their applications, we refer the readers to the articles \cite{ambikasaran2013NlonN} \cite{Khan2022numerical} \cite{Gu2015subspace}. The matrices that possess such hierarchical structures are leveraged to construct various algorithms that can reduce the storage and accelerate the matrix operations \cite{fong2009black}\cite{KHAN2024112786}\cite{lin2011fast} \cite{ambikasaran2019hodlrlib} \cite{ambikasaran2014inverse}. 

In the process of efficiently handling large and dense kernel matrices, hierarchical matrices satisfying weak admissibility criteria (from now on denoted as $ \mathcal{H}_{*} $) have become an essential tool. One of the challenges in dealing with $\mathcal{H}_{*} $ matrices is that the rank of matrices corresponding to neighboring interactions of source and target domains is dependent on the number of source and target particles.


Despite the success of such hierarchical algorithms, most of the theoretical works while studying the rank due to interactions with the neighbors (HODLRdD, HSS2D, etc. \cite{KHAN2024112786} \cite{hackbusch2004hierarchical} \cite{kandappan2023hodlr2d} \cite{borm2003introduction}), assume that the particles (or sources) are placed on a \textbf{uniform grid or quasi-uniform grid}. However, in most practical applications, such an assumption is not often true, as particles rarely align in such structured (\textbf{uniform/quasi-uniform}) patterns. Real-world data, whether coming from various physical simulations, machine learning, or data analysis tasks, tends to follow distributions that are arbitrary in nature, with no consistent pattern. This raises a concern about the applicability and robustness of hierarchical low-rank methods under the arbitrary distribution of particles.

To encounter and model the arbitrary nature of particles, in this article, we consider a more realistic case where particles are randomly distributed in their respective domains according to some suitable probability distribution (see \autoref{subsec: Choice of Probability Distribution} for in-depth discussion), leading to the corresponding kernel matrix being a random matrix (see \autoref{subsec: Generation of Random Kernel Matrix}). Specifically, \textbf{we study the expected growth} (stated in \autoref{thm: expected rank growth in dD}) \textbf{of the random rank $\mclr$} (as defined in \autoref{subsec: Estimating Random Rank: The Random Variable R}) and \textbf{analyze how much this rank $\mclr$ deviates from its expected value} (stated in \autoref{thm: variance of rank in dD}) \textbf{for all possible interactions in $d$-dimensions}. These results offer a more general understanding of algorithms in practical settings. To the best of our knowledge, this study of the rank of the kernel matrices, due to the arbitrary nature of inputs, is a novel contribution to the field. 

\subsection*{Existing Work vs. Our Approach:}  
As mentioned earlier, the rank of kernel matrices has been extensively studied in various applications, from PDEs and inverse problems to machine learning.  In this subsection, we review a few key contributions from existing works that have examined the rank of these kernel matrices, particularly those related to our work, and also highlight how our work differentiates itself from previous works.

Hackbusch et al. \cite{borm2003introduction}\cite{hackbusch2004hierarchical} were among the first to study the rank structure of kernel matrices, interpreting methods like Treecode \cite{barnes1986hierarchical} and FMM \cite{greengard1987fast} \cite{carrier1988fast} as low-rank representations of sub-matrices. Their initial works focused on kernel sub-matrices arising from interactions between well-separated domains, satisfying the standard (or strong) admissibility criterion. One key takeaway from those works is that such \textbf{well-separated interactions naturally lead to low-rank structures} in the corresponding kernel matrix.

In the later works, Hackbusch et al. in ~\cite{hackbusch2004hierarchical} introduced the notion of \textbf{weak-admissibility criterion} in the context of one-dimension and studied the rank of kernel matrices due to interaction between the neighbors in one-dimension (i.e., the vertex sharing case) where they found that the rank is of $\mclo{\log n}$, $n$ being the number of particles in the domains. In the article~\cite{xia2021multi}, Xia extended the weak-admissibility criterion in two dimensions and provided a rough idea of the rank of kernel matrices growth due to the neighboring interactions.

While analyzing the complexity of their algorithms, Ho, Greengard \cite{ho2012fast} and Ho, Lexing \cite{ho2016hierarchical} provide heuristic bounds on the rank. They discussed the rank of interactions between two $d$-dimensional hypercubes sharing a  $(d-1)$-dimensional hypersurface, i.e., vertex-sharing in 1D, edge-sharing in 2D, and face-sharing in 3D, and so on.

Recently, Khan et al. in ~\cite{KHAN2024112786} rigorously proved the rank growth of kernel matrices for all possible neighboring interactions in any dimensions. The result is independent of the choice of kernel functions. However, the main drawback of their work is that to derive the results, they assumed that the particles in each domain are arranged on a uniform or quasi-uniform grid, which is generally not the case in practice. In this work, we relax this assumption and consider \textbf{an arbitrary distribution of particles} in the respective domains, and to model the arbitrary distribution of particles, we consider \textbf{particle distribution} to arise from an underlying \textbf{random distribution}.

Building on these insights, our work extends previous works by analyzing the rank of kernel matrices under arbitrary particle distributions using a probabilistic framework. This novel perspective provides a deeper understanding of the rank growth and its variability, as outlined in the key highlights of this article.


\subsection*{Highlights of The Article:}
The following points are the main highlights of this article, which showcase the unique contributions made in this domain.
\begin{itemize}
    \item To study the behavior of the rank of kernel matrices due to the interactions of \textbf{arbitrary} particles in the neighboring clusters, we assume that the particles are \textbf{randomly distributed} in the respective domains. 
    \item  We have introduced the notion of \textbf{random rank $\mclr$} for the kernel matrices with randomly distributed particles in the respective domains with all possible interactions in $d$-dimensions, which provides a rigorous analysis of
        \begin{itemize}
            \item \textbf{The expected growth of $\mathcal{R}$}. \autoref{thm: expected rank growth in dD} provides deeper insights into how the kernel matrices behave under random conditions, which has not been addressed previously.
            \item \textbf{The change in variance of $\mathcal{R}$}. \autoref{thm: variance of rank in dD} provides a clear understanding of how the rank deviates from its expected value. This analysis helps to explain the stability and variability of hierarchical low-rank algorithms in practical settings.
        \end{itemize}
    \item  To the best of our knowledge, this is the \textbf{first comprehensive study on the expected growth and variance growth} of random kernel matrices for different interactions. Our findings provide a fresh perspective on the robustness and efficiency of hierarchical matrix algorithms in practical settings.
\end{itemize}

\subsection*{Outline of the Article:}

The article is organized as follows. In \autoref{sec: Preliminaries}, we introduce the basic terminologies, definitions, and foundational concepts that are used throughout the article. In \autoref{sec: Main Results}, we state the main theorems and mention the kernel functions that are used to verify our theorems. In \autoref{sec: Fundamental Framework and Problem Setup}, we formally define the problem setup, including the choice of random particle distribution, how the random kernel matrix $K$ is generated, and the random variable $\mclr$. \autoref{sec: Detailed Proof of the Results} presents the proof of theorems on the expected growth of the random rank $\mathcal{R}$ and its variance for different interactions in $d$-dimensions. In \autoref{sec: Numerical Resulsts}, we will provide numerical experiments to validate our theoretical results, focusing on one-, two-, and three-dimensional cases for all possible interactions. Finally, in \autoref{sec: Conclusion}, we summarize the key insights, discuss the implications of our results, and suggest potential directions for future research.

\section{Preliminaries} \label{sec: Preliminaries}
In this section, we are going to discuss some definitions, notations, and lemmas that we will use in the article. We also briefly discuss low-rank approximation of certain kernel matrices using polynomial interpolation and lastly, we will discuss some fundamental probability concepts in our context.

\subsection{Some Notations and Definitions:}

The terminologies that we are going to use frequently in this article are given below.
\paragraph{Source Domain:} The compact set $Y\subset\bbr^d$ is said to be the source domain if it contains the source particles in its interior. We will consider that the particles are arbitrarily distributed within the interior of  $Y$. 
\paragraph{Target Domain:} The compact set $X\subset\bbr^d$ is said to be the target domain if it contains the target particles in its interior. Here also, the target particles are arbitrarily distributed within the interior of $X$. 
Further, we also consider that $\text{int}(X)\cap \text{int}(Y) = \varnothing$, where $\text{int}(X)=$ interior of the set $X$.
\paragraph{Kernel Function:} Throughout the article, we will consider $\mathcal{K}:\text{int}(X)\times \text{int}(Y)\to \bbr$ as the kernel function. The kernel function $\mclk$ encodes the strength of interaction between a particle from the source domain and a particle from the target domain. The choice of a kernel depends on the underlying physical model and the nature of the interaction between particles.

\paragraph{Kernel Matrix:}  The matrix $K\in\bbr^{m\times n}$ whose $\bkt{i,j}^{th}$ entry is given by
\begin{equation}\label{equ: matrix entry K_ij = K(x_i,y_j)}
K_{ij} = \mathcal{K}(x_i,y_j);\quad\text{where } x_i\in \text{int}(X),\, y_j\in \text{int}(Y)     
\end{equation}
is called the kernel matrix. 

\begin{mydfn}{\bf (Numerical $\varepsilon$-rank)}{} \label{def: Numerical rank}
    Given any $\varepsilon>0$, the $\varepsilon$-rank of a matrix $A\in \bbc^{m\times n}$ is denoted and defined as \[r_\varepsilon =\max\bkct{k\in\bkct{1,2,\dotsc,\min\bkct{m,n}} : \frac{\sigma_{k}}{\sigma_1}\geq\varepsilon}, \] where $\sigma_i$'s are singular values of $A$ arranged in decreasing order.
\end{mydfn}
\begin{mydfn}{\bf (Numerical max-rank)}{\label{def: Numerical max-rank}}
    Given any $\varepsilon > 0$, the max-rank of a matrix $A\in \bbc^{m\times n}$ is denoted as $p_\varepsilon$ where $p_\varepsilon = \min\bkct{ \operatorname{rank}\bkt{\tilde{K}} : \tilde{K}\in S_\varepsilon^{(\infty)} } $ and $S_\varepsilon^{(\infty)}$ is defined as \[ S_\varepsilon^{(\infty)} = \bkct{ \tilde{K} \in \bbc^{m\times n} : \magn{ K - \tilde{K}}_{\infty^*} < \varepsilon \magn{K}_{\infty^*} }. \]  
\end{mydfn}
In our theorems, we get an upper bound of $p_\varepsilon$. Here, max-norm of a matrix $A\in \bbc^{m\times n}$ is defined as $ \magn{K}_{\infty^*} = \max\limits_{i,j}\abslt{K\bkt{i,j}} $. Note that \autoref{def: Numerical rank} can also be interpreted as follows:

\textit{Given any $\varepsilon > 0$, the $\varepsilon$-rank of a matrix $A\in \bbc^{m\times n}$ is denoted as $r_\varepsilon$ where $r_\varepsilon = \min\bkct{ \operatorname{rank}\bkt{\tilde{K}} : \tilde{K}\in S_\varepsilon^{(2)} } $ and $S_\varepsilon^{(2)}$ is defined as \[ S_\varepsilon^{(2)} = \bkct{ \tilde{K} \in \bbc^{m\times n} : \magn{ K - \tilde{K}}_2 < \varepsilon \magn{K}_2 }. \]}
\begin{myremark}
    Here, $\varepsilon$-rank and max-rank are related through equivalency between $\magn{\cdot}_2$ and $\magn{\cdot}_{\infty^*}$ (see \autoref{app: Relationship between numerical rank and max-rank}).
\end{myremark}

\subsection{Kernel Matrix Approximation: } Here, we discuss how one can approximate the kernel matrices, and there are many ways to do that. But, in our case, it will be done by approximating the kernel function using Lagrange interpolation. For that purpose, we need the following lemmas.
\begin{mylma} \label{lma:chebinterpolation}
    Let a function $f$ be analytic in $[-1,1]$, and it can be analytically continuable to a Bernstein ellipse for some $\rho>1$ where $|f(x)|\leq M$ for some $M$. Then for any $n \in \bbn$, its Chebyshev interpolants $p_n$ satisfy \[ |f - p_n| \leq \frac{4M\rho^{-n}}{\rho -1}. \]
\end{mylma}
This lemma is proved in [\!\cite{10.5555/3384673}, Theorem 8.2]. We now state a generalized version of \autoref{lma:chebinterpolation} for higher dimensional setting, which is given below
\begin{mylma}\label{lma: cheb interpolation is d dimensions}
     Let $f: V = [-1,1]^d\subset \bbr^d\to \bbr$ be analytic and it has an analytic extension to some generalized Bernstein ellipse $\mathscr{B}\bkt{V,\boldsymbol{\rho}}$, where $\boldsymbol{\rho} = \bkt{\rho,\rho,\dotsc,\rho} $ with $\rho>1$. Now, if $\magn{f}_{\infty^*} =\max\limits_{y\in \mathscr{B}\bkt{V,\boldsymbol{\rho}} } \abslt{f(y)}\leq M $, then for any $n \in \bbn$, its interpolating multivariate polynomial $\Tilde{f}$ satisfy \[ \magn{f - \Tilde{f}}_{\infty^*} \leq 4MV_d\frac{\rho^{-p}}{\rho - 1}, \] where $p$ is a predefined constant and $V_d$ is a constant depending on the dimension $d$ and $\rho$.
\end{mylma} 
The \autoref{lma: cheb interpolation is d dimensions} is discussed in \cite{KHAN2024112786}. Also, a detailed discussion on Bernstein ellipse, generalized Bernstein ellipse, and analytic continuation can be found there. The generalized version of \autoref{lma: cheb interpolation is d dimensions} is stated and proved in \cite{glau2019improved}. However, the above two lemmas provide a way to approximate the kernel function $\mclk(x,y)$. Below, we have shown the approximation of kernel function $\mathcal{K}(x,y)$ along $y$.
\[\mathcal{\Tilde{K}}\bkt{x,y} = \sum_{k\in \mathbb I}\mathcal{K}\bkt{x,y^k} L_k\bkt{y}  \] 
where $\mathbb I$ is the index set of the interpolating points, and $L_k$ is the Lagrange basis. Now, using this approximated kernel function $\mathcal{\Tilde{K}}$, we can get the approximated kernel matrix $\Tilde{K}$ whose $(i,j)^{th}$ entry is given by $\Tilde{K}_{ij} = \mathcal{\Tilde{K}}\bkt{x_i,y_j}$. The matrix is then factorized as $\Tilde{K} = UV^T$,  which is given as follows, where $U\in \bbr^{m\times |\mathbb I|}$ and $V\in \bbr^{n\times|\mathbb I|}$.
\begin{center}
\scalebox{0.6}{ 
    \begin{minipage}{\textwidth} 
        \begin{align*}
            \tilde{K} &= 
            \underbrace{\begin{bmatrix}
                \mathcal{K}(x_1,y^1) & \mathcal{K}(x_1,y^2) & \dotsb & \mathcal{K}(x_1,y^{|\mathbb{I}|}) \\
                \mathcal{K}(x_2,y^1) & \mathcal{K}(x_2,y^2) & \dotsb & \mathcal{K}(x_2,y^{|\mathbb{I}|}) \\
                \vdots & \vdots & \ddots & \vdots \\
                \mathcal{K}(x_m,y^1) & \mathcal{K}(x_m,y^2) & \dotsb & \mathcal{K}(x_m,y^{|\mathbb{I}|}) \\
            \end{bmatrix}}_U
            \times
            \overbrace{\begin{bmatrix}
                L_1(y_1) & L_1(y_2) & \dotsb & L_1(y_n) \\
                L_2(y_1) & L_2(y_2) & \dotsb & L_2(y_n) \\
                \vdots & \vdots & \ddots & \vdots \\
                L_{|\mathbb{I}|}(y_1) & L_{|\mathbb{I}|}(y_2) & \dotsb & L_{|\mathbb{I}|}(y_n) \\
            \end{bmatrix}}^{V^T} 
        \end{align*}
    \end{minipage}}
\end{center}

The rank of the approximated kernel matrix $\Tilde{K}$ is nothing but $|\mathbb I|$. Now, the following lemma guarantees the cardinality of the set $\mathbb{I}$ when the source and target domains are separated by some distance, as mentioned in the following lemma.

\begin{mylma} \label{lma:farfield}
    Let $X$ be the source and $Y$ be the target hyper-cube in $d$-dimensions such that they are at least one hyper-cube away, and let $K$ be the corresponding interaction matrix. Then for any given $ \delta >0$, there exists an approximated $\Tilde{K}$ of rank $p_{ \delta} $ such that $\frac{\magn{K-\Tilde{K}}_{\infty^*}}{\magn{K}_{\infty^*}}< \delta $ with $p_{ \delta }\in \mclo{\log^d\bkt{\frac{c}{ \delta }}}$, where $c$ is a kernel dependent constant.
\end{mylma}
The \autoref{lma: cheb interpolation is d dimensions} is used to prove the \autoref{lma:farfield},  which one can find in the article \cite{KHAN2024112786}, where one can also find a detailed discussion about the kernel-dependent constant $c$. Here, we see that the numerical rank $p_{\delta} $ of the approximated matrix depends on the kernel and desired accuracy $ \delta $, not on the size of the matrix.

    

\section{Main Results} \label{sec: Main Results} 
\autoref{thm: expected rank growth in dD} and \autoref{thm: variance of rank in dD} are among the main contributions of this article. \autoref{thm: expected rank growth in dD} explains the expected growth of the random rank of the kernel matrix due to the interaction of random source and target domains (where source and targets are neighbors) in all possible dimensions, and \autoref{thm: variance of rank in dD} guarantees the growth of variance of $\mclr$. The proof of these theorems can be found in \autoref{sec: Detailed Proof of the Results}. Note that the theorems consider all possible positioning of the source and targets, whereas the \autoref{lma:farfield} deals with only those domains that are separated by a distance. Another point to note is that these two theorems are applicable to an extensive range of kernels, which frequently occur in practice. We have provided numerical validation of the proofs in \autoref{sec: Numerical Resulsts} for the kernels $\mclk_1, \mclk_2, \dotsc, \mclk_7$ mentioned in \autoref{tab: kernel function table} at the end of this section.

We now state the theorem on expected rank growth for all possible hypersurface-sharing interactions.
\begin{boxtheorem}[Expected Rank Growth]\label{thm: expected rank growth in dD}
    Let the compact hyper-cubes $X, Y\subset\bbr^d$ be the source and target domains, respectively, that share a $d'$-dimensional hyper-surface, and each domain contains $n$ i.i.d. random particles in its interior, and the underlying probability distribution is the uniform probability distribution. Let $K$ be the corresponding random kernel matrix. Now, for any given small enough tolerance $ \delta >0$, there exists an approximated random kernel matrix $\Tilde{K}$ such that $\frac{\magn{K-\Tilde{K}}_{\infty^*}}{\magn{K}_{\infty^*}}< \delta $ with random rank $\mclr$ such that
    \begin{enumerate}[(i)]
        \item for $d'=0$, i.e. for vertex-sharing interaction $\bbe \bkbt{\mclr} \in \mclo{p\log_{2^d}(n)}$,
        \item for $d'\neq 0$, i.e. for $d'$-dimensional hyper-surface sharing interaction  $\bbe \bkbt{\mclr} \in \mclo{p n^{d'/d}}$,
    \end{enumerate} 
    where $p$ is a constant depending on $ \delta $ and the kernel.
\end{boxtheorem}

Now, similarly, we state the growth of variance of the random rank $\mclr$ of the kernel matrix $K$ due to all possible interactions.

\begin{boxtheorem}[Variance Growth of the Rank]\label{thm: variance of rank in dD}
    Let the compact hyper-cubes $X, Y\subset\bbr^d$ be the source and target domains, respectively, that share a $d'$-dimensional hyper-surface, and each domain contains $n$ i.i.d. random particles in its interior, and the underlying probability distribution is the uniform probability distribution. Let $K$ be the corresponding random kernel matrix. Now, for any given small enough tolerance $ \delta >0$, there exists an approximated random kernel matrix $\Tilde{K}$ such that $\frac{\magn{K-\Tilde{K}}_{\infty^*}}{\magn{K}_{\infty^*}}< \delta $ with random rank $\mclr$ such that
    \begin{enumerate}[(i)]
        \item for $d'=0$, i.e. for vertex-sharing interaction $\var{\mclr} \in \mclo{(\log\log\log n)^2}$.
        \item for $d'\neq 0$, i.e. for $d'$-dimensional hyper-surface sharing interaction $\var{\mclr}\in \mclo{\bkt{n^{d'/d}\log\log\log n}^2}$.
    \end{enumerate} 
\end{boxtheorem}

We have performed numerical experiments to verify the results, and for that, we have used kernel functions given in \autoref{tab: kernel function table}.

\begin{table}[ht]
    \centering
    \renewcommand{\arraystretch}{1.0} 
    \setlength{\tabcolsep}{11pt}      
    \begin{tabular}{|c|p{5cm}|}
        \hline
        \rowcolor{gray!40}
        \textbf{Sl. No.} &\hfill \textbf{Kernel Functions}\hfill\textcolor{gray!40}{.} \\
        \hline
        \hline
        1 & $\mclk_1\bkt{x,y} = \dfrac{1}{\magn{x-y}_2}$ \\
        \hline
        \rowcolor{gray!20}
        2 & $\mclk_2\bkt{x,y} = \log\bkt{\magn{x-y}_2}$ \\
        \hline
        3 & $\mclk_3\bkt{x,y} = \sin\bkt{\magn{x-y}_2}$ \\
        \hline
        \rowcolor{gray!20}
        4 & $\mclk_4\bkt{x,y} = \dfrac{\exp\bkt{i\magn{x-y}_2}}{\magn{x-y}_2}$ \\
        \hline
        5 & $\mclk_5\bkt{x,y} = \dfrac{1}{\sqrt{1+\magn{x-y}_2}}$ \\
        \hline
        \rowcolor{gray!20}
        6 & $\mclk_6\bkt{x,y} = \exp\bkt{-\magn{x-y}_2}$ \\
        \hline
        7 & $\mclk_7\bkt{x,y} = \magn{x-y}_2$ \\
        \hline
    \end{tabular}
    \caption{Table of Kernel Functions where $x\in\text{int}(X)$ and $y\in\text{int}(Y)$. Here $\mclk_1$ is the Green's function for the 3D Laplace operator; $\mclk_2$ is the Green's function for the 2D Laplace operator; $\mclk_4$  is the 3D Helmholtz kernel with wave number 1; $\mclk_6$ is the Matérn covariance kernel, and $\mclk_7$ is the poly-harmonic radial basis function. }
    \label{tab: kernel function table}
\end{table}

\section{Fundamental Framework and Problem Setup} \label{sec: Fundamental Framework and Problem Setup}
In this section, we provide a foundational discussion of the probability distribution and key concepts that are essential for understanding the results presented in this article. The distribution and the random variables defined here play a crucial role in the proofs of the theorems stated in \autoref{sec: Main Results} and provide the necessary framework for analyzing the behavior of the random variable $\mclr$ of our interest.

For the sake of simplicity, from now onward, we are going to consider the hyper-cube $Y = [0, l]^d$ as the source domain and the hyper-cube $X = [0,l]^{d'}\times [-l, 0]^{d-d'}$ as the target domain in the $d$-dimensions. Here, for different values of $d'$, the source and target will encounter different interactions, which are given below
\begin{itemize}
    \item For $d'=-1$, the source $Y$ and the target $X$ are \textbf{far-field}\footnote{For the sake of completeness, we use the notation $d'=-1$ to represent far-field domains. There is nothing special to choose $-1$ here.} domains, where far-field is defined as $X$ and $Y$ are at least one hyper-cube away, i.e. $\text{dist}\bkt{X,Y}\geq \eta\min\bkct{\text{diameter}(X),\text{diameter}(Y)} $, for some $\eta>0$. The far-field domains in one, two, and three dimensions are shown in the following \autoref{fig: far-field in d-D}.
    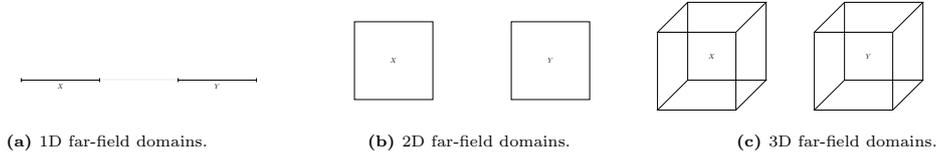
\begin{figure}[H]
        \centering
        \resizebox{0.7\textwidth}{!}{
        \begin{tikzpicture}
            \begin{scope}
                \draw (-4,0) -- (8,0);
                \draw[white, line width=1pt] (0,0) -- (4,0);
                \draw (-4,0.1) -- (-4,-0.1);
                \draw (0,0.1) -- (0,-0.1);
                \draw (4,0.1) -- (4,-0.1);
                \draw (8,0.1) -- (8,-0.1);
                \node at (-2,-0.3) {$X$};
                \node at (6,-0.3) {$Y$};
            \end{scope}
            \begin{scope}[shift ={(13,-1)}]
                \draw (0,0) rectangle (4,4);
                \node at (2,2) {$X$};
                \draw (8,0) rectangle (12,4);
                \node at (10,2) {$Y$};
            \end{scope}
            \begin{scope}[shift ={(30,0)}]
                \def\lll{4}
                \drawcube{0}{0}{0}{\lll}  
                \drawcube{8}{0}{0}{\lll}    
                \node at (2,2,2) {$X$};
                \node at (10,2,2) {$Y$};
            \end{scope}
        \end{tikzpicture}}
        \resizebox{0.85\textwidth}{!}{
        \begin{minipage}{0.32\textwidth}
        \centering
        \subcaption{1D far-field domains.}
        \label{fig: 1d far-fields plain}
        \end{minipage}
        \hfill
        \begin{minipage}{0.32\textwidth}
            \centering
            \subcaption{2D far-field domains.}
            \label{fig: 2d far-field plain}
        \end{minipage}
        \hfill
        \begin{minipage}{0.33\textwidth}
            \centering
            \subcaption{3D far-field domains.}
            \label{fig: 3d far-field plain}
        \end{minipage}}
        \caption{1D, 2D, 3D far-field domains}
        \label{fig: far-field in d-D}
    \end{figure}
    \item For $d'$ = 0, the source $Y$ and the target $X$ shares a vertex. \autoref{fig: vertex sharing in d-D} illustrates the vertex-sharing domains in one, two, and three dimensions.
    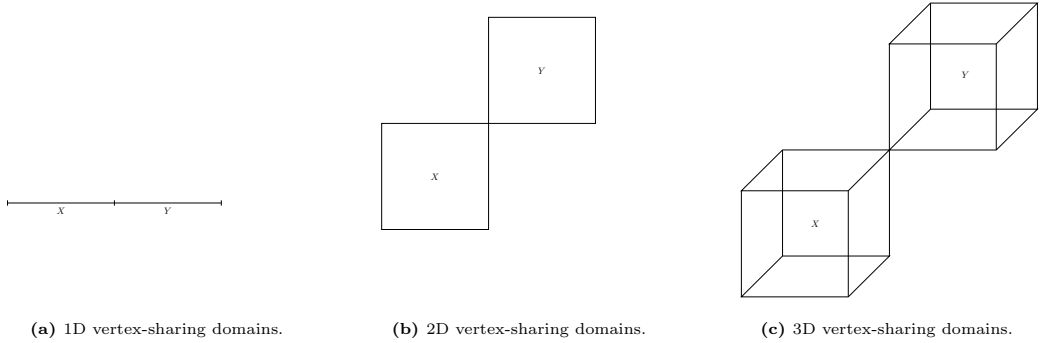
\begin{figure}[H]
        \centering
        \resizebox{0.8\textwidth}{!}{
        \begin{tikzpicture}
            \begin{scope}
                \draw (-4,0) -- (4,0);
                \draw (-4,0.1) -- (-4,-0.1);
                \draw (0,0.1) -- (0,-0.1);
                \draw (4,0.1) -- (4,-0.1);
                \node at (-2,-0.3) {$X$};
                \node at (2,-0.3) {$Y$};
            \end{scope}
            \begin{scope}[shift ={(10,-1)}]
                \draw (0,0) rectangle (4,4);
                \node at (2,2) {$X$};
                \draw (8,8) rectangle (4,4);
                \node at (6,6) {$Y$};
            \end{scope}
            \begin{scope}[shift ={(25,-2)}]
                \def\lll{4}
                \drawcube{0}{0}{0}{\lll}  
                \drawcube{6.5}{6.5}{2.5}{\lll}    
                \node at (2,2,2) {$X$};
                \node at (8,8,3) {$Y$};
            \end{scope}
        \end{tikzpicture}}
        \resizebox{0.85\textwidth}{!}{
        \begin{minipage}{0.32\textwidth}
        \centering
        \subcaption{1D vertex-sharing domains.}
        \label{fig: 1d vertex-sharing plain}
        \end{minipage}
        \hfill
        \begin{minipage}{0.32\textwidth}
            \centering
            \subcaption{2D vertex-sharing domains.}
            \label{fig: 2d vertex-sharing plain}
        \end{minipage}
        \hfill
        \begin{minipage}{0.33\textwidth}
            \centering
            \subcaption{3D vertex-sharing domains.}
            \label{fig: 3d vertex-sharing plain}
        \end{minipage}}
        \caption{1D, 2D, 3D vertex-sharing domains}
        \label{fig: vertex sharing in d-D}
    \end{figure}
    \item For $d'$ = 1, the source $Y$ and the target $X$ shares an edge. This case will arise when the dimension of the domains is at least two. In two and three dimensions, this case is illustrated in \autoref{fig: edge sharing in d-D}.
       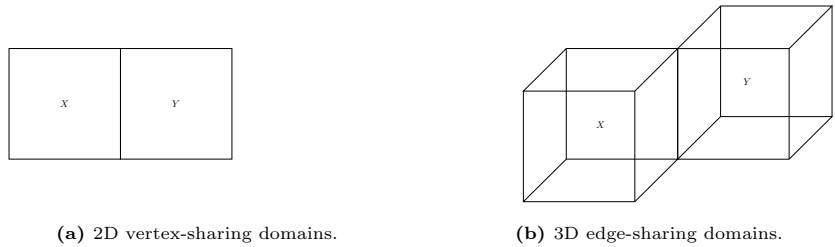
\begin{figure}[H]
           \centering
           \resizebox{0.64\textwidth}{!}{
            \begin{tikzpicture}
                \begin{scope}
                    \draw (0,0) rectangle (4,4);
                    \node at (2,2) {$X$};
                    \draw (4,0) rectangle (8,4);
                    \node at (6,2) {$Y$};
                \end{scope}
                \begin{scope}[shift ={(20,0)}]
                    \def\lll{4}
                    \drawcube{0}{0}{0}{\lll}  
                    \drawcube{6.5}{2.5}{2.5}{\lll}    
                    \node at (2,2,2) {$X$};
                    \node at (7,3.3,1.3) {$Y$};
                \end{scope}
            \end{tikzpicture}}
            \resizebox{0.7\textwidth}{!}{
            \begin{minipage}{0.36\textwidth}
                \centering
                \subcaption{2D vertex-sharing domains.}
                \label{fig: 2d edge-sharing plain}
            \end{minipage}
            \hfill
            \begin{minipage}{0.36\textwidth}
                \centering
                \subcaption{3D edge-sharing domains.}
                \label{fig: 3d edge-sharing plain}
            \end{minipage}}
           \caption{2D, 3D edge-sharing domains}
           \label{fig: edge sharing in d-D}
       \end{figure}
    \item For $d'$ = 2, the source $Y$ and the target $X$ shares a face. This case will arise when the dimension of the domains is at least three. The face-sharing domain in three dimensions is shown below in \autoref{fig: face sharing in 3D}.

    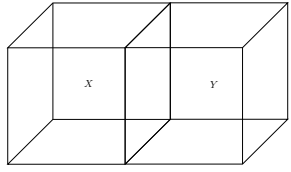
\begin{figure}[H]
        \centering
        \resizebox{0.22\textwidth}{!}{
        \begin{tikzpicture}
            \def\lll{4}
                \drawcube{0}{0}{0}{\lll}  
                \drawcube{6.5}{2.5}{6.48}{\lll}    
                \node at (2,2,2) {$X$};
                \node at (6,1.7,1.3) {$Y$};
        \end{tikzpicture}}
        \caption{3D face-sharing domains}
        \label{fig: face sharing in 3D}
    \end{figure}

    \item For $d' \geq 3$, the source $Y$ and the target $X$ shares a $d'$-dimensional hyper-surface. 
\end{itemize}
\subsection{Choice of Probability Distribution} \label{subsec: Choice of Probability Distribution}

We begin this subsection by introducing the probability distribution that forms the foundation of our analysis as we model the `arbitrary distribution of particles' with a proper probability distribution. In each of domains $X$ and $Y$, $n$ independent and identically distributed (i.i.d) particles are drawn from \textbf{Uniform Probability distributions} $\mathcal{U}_{X}$ and $\mathcal{U}_{Y}$ respectively. The choice of Uniform Probability Distribution is driven by its suitability and simplicity for modeling the random matrix structure considered in this study.
Apart from simplicity in deriving analytical results, several factors influenced this choice, which are summarized as follows:
\begin{enumerate}[(i)]
    \item  Instead of being tied to a fixed grid, uniform probability distribution allows for flexible grid placement during the discretization of the domains.
    \item Uniform distribution ensures that every possible particle within the domain is equally likely to be chosen.
    \item The configuration that the particles sampled randomly from the uniform probability distribution can be interpreted as a ``perturbed" version of a uniform grid, which aligns with the principles of \textbf{Smoothed Analysis} \cite{spielman2009smoothed} in the following ways:
    \begin{enumerate}
        \item Choosing uniform probability distribution, we can implicitly introduce one kind of random perturbation to an ideal grid\footnote{A detailed discussion on `Uniform Distribution as Randomly Uniform Perturbed Grid' can be found in \autoref{app: Uniform Distribution as Randomly Uniform Perturbed Grid}}, as each particle selected from the domain could vary slightly from an exact grid particle. This randomness in particle selection in the domains aligns with the smooth perturbations considered in the smooth analysis.
        \item The use of a uniform probability distribution can provide more realistic estimates of algorithm performance, as it takes a range of possible input particles into account rather than specific, deterministic ones.
        \item By choosing uniform probability distribution, we are spreading the likelihood evenly across the domains, making it less likely to encounter those extreme cases that result in an analysis that better reflects typical performance (which is quite similar to average-case analysis of inputs) rather than focusing on the rare, worst-case analysis of inputs particles.    
    \end{enumerate}
\end{enumerate}

\subsection{Generation of Random Kernel Matrix} \label{subsec: Generation of Random Kernel Matrix}
Following the discussion on the probability distribution, we now turn to the generation of the random kernel matrix $K$. As particles are uniformly distributed in their respective domains and the kernel matrix reflects the interactions between these random particles through a kernel function $\mclk$, forms the random entries of the matrix $K$ as discussed in \autoref{equ: matrix entry K_ij = K(x_i,y_j)}. 
Although the kernel matrix $K$ exhibits randomness, this randomness is structured\footnote{{\bf Structured Randomness} involves randomness that is governed by some underlying structure, pattern, or correlation between the elements. This structure can result from dependencies, constraints, or specific relationships defined by a process or function.} and is not entirely arbitrary. This is because:
\begin{enumerate}[i.]
    \item The entries $K_{ij}$ of the kernel matrix $K$ are random variables whose distribution depends not only on $\mathcal{U}_X$ and $\mathcal{U}_Y$ but also on the kernel function $\mclk$.
    \item The distribution and correlation structure of the entries $K_{ij}$ of the kernel matrix $K$ are strongly influenced by the geometrical configurations of the domains $X$ and $Y$. 
    \item The distribution of $K_{ij}$ will often not be uniform (often highly structured). For example, Suppose we choose the kernel function $\mclk_6$. In that case, the entries of $K$ will be related to the exponential of the Euclidean distance between the random particles, leading to a non-uniform distribution.
\end{enumerate}

\subsection{Rank as a Random Variable}
Once the random kernel matrix $K$ is generated, one of the key aspects of its structure is its rank. Unlike deterministic kernel matrices, where the rank is fixed based on the matrix dimensions and locations of the source and target domains, the rank of the kernel matrix derived from random source and target particles is a random variable. Several factors influence this random variable:
\begin{enumerate}[i.]
    \item  The size of the kernel matrix $K$ (i.e., the number of random source and target particles)
    \item The spatial configuration of source and target domains has a substantial impact on the rank of $K$. The rank can be significantly different in cases when the domains are far-field or share a $d'$-dimensional hyper-surface.
    \item The rank of the kernel matrix is not only depends on the size of $K$ and geometry of the source and target domains but also heavily influenced by the choice of the kernel function $\mclk$. For example, for the choice kernel function of $\mclk_3$ from \autoref{tab: kernel function table}, regardless of the geometry or distribution of the source and target particles in one dimension, the rank of the kernel matrix is always 2.
\end{enumerate}    
The above-outlined points can be verified in the numerical results presented in \autoref{sec: Numerical Resulsts}. As determining the exact distribution of the rank is quite challenging, we introduce the random variable $\mclr$, (given in \autoref{equ: random rank in dD}), which serves as an upper bound of the rank of kernel matrix under random inputs. Throughout this article, we investigate the behavior of $\mclr$ to understand how the rank of the kernel matrix varies and obtain bounds on the first couple of moments of the random variable.

\subsection{Basic Probability Theory}\label{subsec: Basic Probability Theory} 
In this subsection, we explore fundamental probability concepts that are crucial for our analysis, and along with this, we define some required tools.

\subsubsection{Distribution of particles in a domain and its sub-domains.} Let $n$ number of independent and identically distributed (i.i.d) particles $x_1,x_2,\dotsc,x_n$ fall in the hyper-cube $V = [a,b]^d$ under 
the uniform probability distribution\footnote{In \autoref{app: Generalization to Arbitrary Probability Distributions}, the consideration of more general probability distributions is discussed.}. The reason behind choosing this particular distribution is discussed in \autoref{subsec: Choice of Probability Distribution}.

We define the random variable $N\bkt{V'}$ as the number of particles that fall within a sub-hyper-cube $V' = [c_1,d_1]^d$ of $V$, where $[c_1,d_1] \subseteq [a,b]$. Then, the probability of having exactly $k$ particles in $V'$ is given by 
\begin{equation}\label{equ: binomial distribution pmf}
    \bbpb{N\bkt{V'}=k} = \binom{n}{k}\mathtt{q}^k\bkt{1-\mathtt{q}}^{n-k}; \quad \text{where }\mathtt{q} = \frac{\bkt{d_1-c_1}^d}{\bkt{b-a}^d}.
\end{equation}

We now define another random variable that we are going to use frequently in the article. The random variable is defined as follows.
\begin{equation}\label{equ: minNkp}
     \mclz_n^{V'} = \min\bkct{N\bkt{V'},p}, \text{where $p$ is a constant}. 
\end{equation}
Here, $\mclz_n^{V'}$ takes any natural number in between 0 and $p$ and for any $i\in \bkct{0,1,\dotsc,p}$, we have 
\begin{equation} \label{equ: minNkp pfm}
    \bbpb{\mclz_n^{V'}=i}=\begin{cases}
            \bbp\bkt{N\bkt{V'}=i} &\text{if } i<p\\
            \bbp\bkt{N\bkt{V'}\geq p} &\text{if } i=p\\
        \end{cases}
\end{equation}

We now demonstrate an application of the trinomial distribution in our setting as follows: Let us consider two non-intersecting sub-hyper-cubes $V'=[c_1,d_1]^d$ and $V''=[c_2,d_2]^d$ of $V$. Then the probability that exactly $l,m$ particles fall in the respective sub-hyper-cubes is given by
\begin{equation}\label{equ: trinomial distribution pmf}
    \bbpb{N\bkt{V'}=l, N\bkt{V''}=m} =\binom{n}{l}\binom{n-l}{m}\mathtt{q_1}^l\mathtt{q_2}^m(1 - \mathtt{q_1}- \mathtt{q_2})^{(n-l-m)};\text{where }\mathtt{q}_i = \frac{\bkt{d_i-c_i}^d}{\bkt{b-a}^d}, i=1,2.
\end{equation}
Now, one more application of trinomial distribution can be given for the random variable defined in \autoref{equ: minNkp} as follows:
\begin{equation}\label{equ: minZniZnj pdf}
    \bbpb{\mclz_n^{V'}=l,\mclz_n^{V''}=m}=\begin{cases}
        \bbpb{N\bkt{V'}=l,N\bkt{V''}=m} &\text{if } l<p,m<p\\
        \bbpb{N\bkt{V'}\geq p, N\bkt{V''}=m} &\text{if } l=p,m<p\\
        \bbpb{N\bkt{V'} =l, N\bkt{V''}\geq p} &\text{if } l<p,m=p\\
        \bbpb{N\bkt{V'} \geq p, N\bkt{V''}\geq p} &\text{if } l=p,m=p\\
    \end{cases}
\end{equation}



\subsubsection{Moments and Dependencies of Random Variables.}
We now discuss the expectation and the variance of the random variables that we have defined in \autoref{equ: binomial distribution pmf} and \autoref{equ: minNkp}. As, $N\bkt{V'}$ is nothing but a binomial distribution with parameters $n, \mathtt{q}$, corresponding expectation and variance is 
\begin{equation}\label{equ: expectation of N(a,b)}
    \bbeb{N\bkt{V'}} = n\mathtt{q}\quad\quad\text{ and }\quad\quad \var{N\bkt{V'}}= n\mathtt{q}(1-\mathtt{q}).
\end{equation}
Now, the expectation of the random variable $\mclz_n^{V'}$ is given below
\begin{align} \label{equ: expectation of Znab}
    \bbeb{\mclz_n^{V'}} = \sum_{i=0}^p i\,\bbpb{\mclz_n^{V'}=i} 
                           &= p + \sum_{i=0}^{p-1} \bkt{i-p}\bbpb{\mclz_n^{V'}=i}\notag\\
                           &= p + \sum_{i=0}^{p-1} \bkt{i-p}\bbpb{N\bkt{V'}=i}\\
                           &= p - \sum_{i=0}^{p} i\,\bbpb{N\bkt{V'}=i} \label{equ: expectation of Znab compact form}
\end{align}
and the variance can be derived as follows
\begin{align}\label{equ: variance Znab}
    \var{\mclz_n^{V'}} = \var{p-\mclz_n^{V'}}
                          &=\sum_{i=0}^{p} \bkt{p-i}^2\,\bbpb{\mclz_n^{V'}=i} - \bkt{\sum_{i=0}^{p} \bkt{p-i}\,\bbpb{\mclz_n^{V'}=i}}^2 \notag \\
                          &= \sum_{i=0}^p i^2\,\bbpb{N\bkt{V'}=i} - \bkt{\sum_{i=0}^p i\,\bbpb{N\bkt{V'}=i}}^2
\end{align}
Now, for any non-intersecting sub-hyper-cubes $V', V''$ of $V$, the expectation of the product and covariance of the random variable defined in \autoref{equ: binomial distribution pmf}, can be easily derived and is given below \begin{equation}
    \bbeb{N\bkt{V'}N\bkt{V''}} = n(n-1)\mathtt{q_1q_2}\quad \text{ and } \quad \cov{N\bkt{V'},N\bkt{V''}} = -n\mathtt{q_1q_2}
\end{equation} 
\begin{myremark}
    The negative covariance indicates that if the number of particles increased in one of the sub-hyper-cubes, then the number of particles should decrease in the other sub-hyper-cube.  
\end{myremark}
Now, the covariance between $\mclz_n^{V'},\mclz_n^{V''}$ for the non intersecting sub-hyper-cubes $V', V''$ of $V$, can be derived as follows  \begin{align} \label{equ: cov(ZniZnj)}
    \cov{\mclz_n^{V'},\mclz_n^{V''}} &=\cov{p-\mclz_n^{V'},p-\mclz_n^{V''}}\notag\\  
                              &= \sum_{l=0}^{p}\sum_{m=0}^{p} lm\bkt{\bbpb{N\bkt{V'}=l,N\bkt{V''}=m}-\bbpb{N\bkt{V'}=l}\bbpb{N\bkt{V''}=m}}
\end{align}
Similarly, the covariance between $\mclz_n^{V'}$ and $N\bkt{V''}$ for the non intersecting sub-hyper-cubes $V', V''$ of $V$, can be derived as {\small\begin{equation}\label{equ: cov(Znk_Mkappa)}
    \cov{\mclz_n^{V'},N\bkt{V''}} = \sum_{l=0}^{p}\sum_{m=0}^{n} lm\bkt{\bbpb{N\bkt{V'}=l,N\bkt{V''}=m}-\bbpb{N\bkt{V'}=l}\bbpb{N\bkt{V''}=m}}
\end{equation}}

\subsubsection{Continuous Approximation of Discrete Distributions. }
In many practical scenarios of discrete random variables, particularly while dealing with large sample sizes or simplifying the computation of probabilities, we employ the method of continuous approximation. The \textit{Central Limit Theorem} \cite{ross2009probability} \cite{Bhattacharya1968Berry} is the primary tool that ensures that we can do so. We now briefly discuss the normal approximation to the binomial distributions.
\begin{mylma}
    Let $S_n$ be a binomial random variable with parameter $n,\mathtt{q}$, i.e. $S_n\sim \textbf{Binomial}(n,\mathtt{q})$, then for $a,b\in\bbr$ \[ \bbpb{a\leq \frac{S_n - n\mathtt{q}}{\sqrt{n\mathtt{q}(1-\mathtt{q})}} \leq b }\to \bbpb{a\leq Z \leq b}, \] as $n\to\infty$, where $Z$ is the standard normal distribution, i.e. $Z\sim\mcln\bkt{0,1}$.
\end{mylma} 
\noindent To ensure better approximation, we need the following assumptions
\begin{enumerate}[(i)]
    \item the quantities $n\mathtt{q}$ and $n(1-\mathtt{q})$ should be a large value (some authors suggest that if both the values are at least 10, then we can get a good approximation, and few suggest 5 is sufficient).
    \item continuity correction is required for best approximation\cite{ross2009probability}.
\end{enumerate}
Now, using Berry-Esseen theorem \cite{Feller1971AnIntro}, we can approximate the probability defined in \autoref{equ: binomial distribution pmf} using the CDF of $\mcln(0,1)$, i.e. $\displaystyle\Phi\bkt{x}=\frac{1}{\sqrt{2\pi}}\int_{-\infty}^x e^{-t^2/2} dt$ and the corresponding error of approximation is \[ \abslt{\bbpb{N(V')=k} - \int_{a_{k,\mathtt{q}}^{(n)}}^{b_{k,\mathtt{q}}^{(n)}} \frac{1}{\sqrt{2\pi}} e^{-t^2/2} dt  } \leq \frac{ 2C (1-2\mathtt{q} + 2\mathtt{q}^2 )}{\sqrt{n \mathtt{q} (1-\mathtt{q})}} \in \mclo{\frac{1}{\sqrt{n}}}, \] where $ a_{k,\mathtt{q}}^{(n)} = \frac{k-0.5-n\mathtt{q}}{\sqrt{n\mathtt{q}\bkt{1-\mathtt{q}}}} $ and $ b_{k,\mathtt{q}}^{(n)} = \frac{k+0.5-n\mathtt{q}}{\sqrt{n\mathtt{q}\bkt{1-\mathtt{q}}}} $. Also note that $C$ is independent of all quantities in the above expression. \\

Now similarly, using the Berry-Esseen theorem for multivariate cases \cite{Bentkuus2004ALyap}, we can approximate the trinomial distribution defined in \autoref{equ: trinomial distribution pmf}, using the bivariate normal distribution \cite{teo2024measuring} and the corresponding error of approximation is given by \[ \abslt{ \bbpb{N\bkt{V'}=l, N\bkt{V''}=m} - \iint\limits_{\mathcal{D}} f(x,y)dx dy } \leq \frac{4M}{\sqrt{n}}, \quad \text{ where $M$ is some constant. } \]
where 
\begin{enumerate}
    \item $l+m\leq n$.
    \item $f(x,y) = \dfrac{1}{2\pi\sqrt{1-\rho^2}}e^{-Q\bkt{x,y}/2\bkt{1-\rho^2}}$ is the pdf of bivariate normal distribution with mean vector $\mu= \begin{bmatrix}0\\0\end{bmatrix}$ and covariance matrix $\Sigma= \begin{bmatrix}1 &\rho\\\rho&1\end{bmatrix}$,  where
    \begin{enumerate}
        \item $Q(x,y)= x^2+y^2-2\rho xy$.
        \item $\rho= -\sqrt{\dfrac{\mathtt{q_1}\mathtt{q_2}}{\bkt{1 - \mathtt{q_1}}\bkt{1 - \mathtt{q_2}}}}$ be the correlation coefficient between $N\bkt{V'}$, $N\bkt{V''}$.
    \end{enumerate}
    \item $\mathcal{D}$ is a rectangular region in $xy$ plane containing the point $\bkt{\frac{l-n\mathtt{q_1}}{\sqrt{n\mathtt{q_1}\bkt{1-\mathtt{q_1}}}}, \frac{m-n\mathtt{q_2}}{\sqrt{n\mathtt{q_2}\bkt{1-\mathtt{q_2}}}}}$ defined as \[\mathcal{D} = \bkct{ (x,y) : \frac{l -0.5-n\mathtt{q_1}}{\sqrt{n\mathtt{q_1}\bkt{1-\mathtt{q_1}}}}\leq x \leq \frac{l + 0.5 -n\mathtt{q_1}}{\sqrt{n\mathtt{q_1}\bkt{1-\mathtt{q_1}}}}; \, \frac{m -0.5-n\mathtt{q_2}}{\sqrt{n\mathtt{q_2}\bkt{1-\mathtt{q_2}}}}\leq y \leq \frac{m + 0.5 -n\mathtt{q_2}}{\sqrt{n\mathtt{q_2}\bkt{1-\mathtt{q_2}}}} } \]
\end{enumerate}

\subsection{Hierarchical Subdivision of the Source Domain}
To define the random variable $\mclr$, we use hierarchical subdivision on the source domain $Y$. This subdivision is done in a way that adapts to the relative positioning of the source and target domains. To get a clear understanding of this, we discuss the one-dimensional case first, and then we move to the general $d$-dimensional case.

\subsubsection*{Subdivision in One Dimension} 
The hierarchical subdivision will take place in this dimension only when the source $Y$ and the target $X$ share a vertex. The source domain $Y$ is hierarchically subdivided using an adaptive binary tree as shown in \autoref{fig: hierarchical subdibvision}, where each level of subdivision halves the domain:
\begin{itemize}
    \item At \textbf{level 1}, $Y$ is subdivided into two equal parts $Y_1$, which shares a vertex with the $X$ and $Y_{1,1}$, which does not.
    \item At \textbf{level 2}, $Y_1$ is now subdivided again into two equal parts $Y_2$ and $Y_{2,1}$ where $Y_2$ shares a vertex with the $X$, and so on. 
\end{itemize}
The hierarchical subdivision continues up to level $\kappa$, where $\kappa = \left\lfloor\log_2 n\right\rfloor$. The process can be summarized as: \begin{align*}
        Y = \underbrace{Y_1\cup Y_{1,1}}_{\text{Level 1}} = \overbrace{ Y_2 \bigcup_{i=1}^2 Y_{i,1}}^{\text{Level 2}} = \dotsb = \overbrace{ Y_\kappa\bigcup_{i=1}^\kappa Y_{i,1}}^{\text{Level }\kappa}.
    \end{align*}
        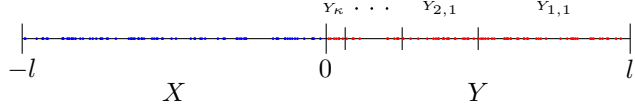
\begin{figure}[htbp]
            \centering
            \begin{tikzpicture}
                \draw (0,0) -- (8,0);
                \draw (0,-0.2) -- (0,0.2);
                \draw (6,-0.15) -- (6,0.15);
                \draw (5,-0.15) -- (5,0.15);
                \draw (4.25,-0.15) -- (4.25,0.15);
                \draw (4,-0.2) -- (4,0.2);
                \draw (8,-0.2) -- (8,0.2);

                \foreach \i in {1,...,80}{
                    \pgfmathsetmacro{\randomx}{6+2*rand}
                    \fill[red] (\randomx,0) circle(0.5pt);
                }
                \foreach \i in {1,...,80}{
                    \pgfmathsetmacro{\randomx}{6+2*rand}
                    \fill[blue] (\randomx-4,0) circle(0.5pt);
                }
                
                \fill (4.4, 0.4) circle (0.5pt);
                \fill (4.6, 0.4) circle (0.5pt);
                \fill (4.8, 0.4) circle (0.5pt);
            
                \node at (7,0.4) {\tiny$Y_{1,1}$};
                \node at (5.5,0.4) {\tiny$Y_{2,1}$};
                \node at (4.1,0.4) {\tiny$Y_{\kappa}$};
                \node at (2,-0.7) {$X$};
                \node at (6,-0.7) {$Y$};
                \node at (0,-0.4) {$-l$};
                \node at (3.99,-0.4) {$0$};
                \node at (8,-0.4) {$l$};
            \end{tikzpicture}
            \caption{Sub-division of the source $Y=[a,b]$ up to the level $\kappa$.}
            \label{fig: hierarchical subdibvision}
        \end{figure}

\subsubsection*{General Case: Subdivision in $d$ Dimensions}
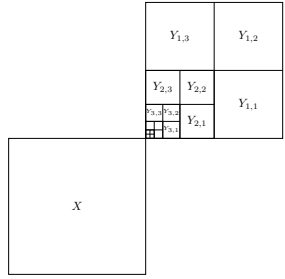
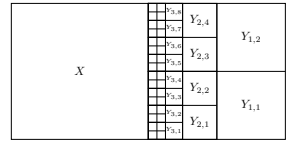
\begin{figure}[H]
    \centering
    \resizebox{0.6\textwidth}{!}{
    \begin{tikzpicture}
        \begin{scope}
            \draw (0,0) rectangle (4,4);
            \node at (2,2) {$X$};
            
            \drawsubdividedsquare{4}{4}{4}  
            \node at (7,5) {$Y_{1,1}$};
            \node at (7,7) {$Y_{1,2}$};
            \node at (5,7) {$Y_{1,3}$};
            
            \drawsubdividedsquare{4}{4}{2}
             \node at (5.5,4.5) {$Y_{2,1}$};
            \node at (5.5,5.5) {$Y_{2,2}$};
            \node at (4.5,5.5) {$Y_{2,3}$};
            
            \drawsubdividedsquare{4}{4}{1}
            \node at (4.75,4.25) {\tiny $Y_{3,1}$};
            \node at (4.75,4.75) {\tiny $Y_{3,2}$};
            \node at (4.25,4.75) {\tiny $Y_{3,3}$};
            
            \drawsubdividedsquare{4}{4}{0.5}
            \drawsubdividedsquare{4}{4}{0.25}    
        \end{scope}
        \begin{scope}[shift={(14.5,1)}]
            \draw (0,0) rectangle (4,4);
            \node at (2,2) {$X$};
            
            \drawsubdividedsquare{4}{0}{4}
            \node at (7,1) {$Y_{1,1}$};
            \node at (7,3) {$Y_{1,2}$};
            
            \drawsubdividedsquare{4}{0}{2}
            \drawsubdividedsquare{4}{2}{2}
            \node at (5.5,0.5) {$Y_{2,1}$};
            \node at (5.5,1.5) {$Y_{2,2}$};
            \node at (5.5,2+0.5) {$Y_{2,3}$};
            \node at (5.5,2+1.5) {$Y_{2,4}$};

            \drawsubdividedsquare{4}{0}{1}
            \drawsubdividedsquare{4}{1}{1}
            \drawsubdividedsquare{4}{2}{1}
            \drawsubdividedsquare{4}{3}{1}
            \node at (4.75,0.25) {\tiny $Y_{3,1}$};
            \node at (4.75,3*0.25) {\tiny $Y_{3,2}$};
            \node at (4.75,5*0.25) {\tiny $Y_{3,3}$};
            \node at (4.75,7*0.25) {\tiny $Y_{3,4}$};
            \node at (4.75,9*0.25) {\tiny $Y_{3,5}$};
            \node at (4.75,11*0.25) {\tiny $Y_{3,6}$};
            \node at (4.75,13*0.25) {\tiny $Y_{3,7}$};
            \node at (4.75,15*0.25) {\tiny $Y_{3,8}$};

            \drawsubdividedsquare{4}{0}{0.5}
            \drawsubdividedsquare{4}{0.5}{0.5}
            \drawsubdividedsquare{4}{1}{0.5}
            \drawsubdividedsquare{4}{1.5}{0.5}
            \drawsubdividedsquare{4}{2}{0.5}
            \drawsubdividedsquare{4}{2.5}{0.5}
            \drawsubdividedsquare{4}{3}{0.5}
            \drawsubdividedsquare{4}{3.5}{0.5}
        \end{scope}
        \end{tikzpicture}}
    \begin{minipage}{0.48\textwidth}
        \centering
        \subcaption{Vertex-sharing interaction.}
        \label{fig: 2d vertex-sharing}
    \end{minipage}
    \hfill
    \begin{minipage}{0.48\textwidth}
        \centering
        \subcaption{Edge-sharing interaction.}
        \label{fig: 2d edge-sharing}
    \end{minipage}

    \caption{Sub-division of the source $Y$ in two dimensions.}
    \label{fig: 2d combined panels}
        
\end{figure}

\begin{figure}[H]
    \centering 
    \resizebox{0.65\textwidth}{!}{
    \begin{tikzpicture}
    \begin{scope}
        \def\lll{4}
        \drawcube{0}{0}{0}{\lll}  
        \drawsubcubes{6.5}{6.5}{2.5}{4}
        \drawsubcubes{5.81}{5.81}{2.675}{2} 
        \drawsubcubes{5.91}{5.92}{4}{1}
    \end{scope}
    \begin{scope}[shift={(20,0.5)}]
        \def\lll{4}
        \drawcube{0}{0}{0}{\lll}  
        \drawsubcubes{6.5}{6.5-4}{2.5}{4}
        \drawsubcubes{5.81}{5.81-4}{2.675}{2} 
        \drawsubcubes{5.81}{5.81-2}{2.675}{2} 
    \end{scope}    
    \end{tikzpicture}}
    \begin{minipage}{0.48\textwidth}
        \centering
        \subcaption{Vertex-sharing interaction.}
        \label{fig:subfig1}
    \end{minipage}
    \hfill
    \begin{minipage}{0.48\textwidth}
        \centering
        \subcaption{Edge-sharing interaction.}
        \label{fig:subfig2}
    \end{minipage}

    \resizebox{0.25\textwidth}{!}{
    \begin{tikzpicture}
    \begin{scope}
        \def\lll{4}
        \drawcube{0}{0}{0}{\lll}  
        \drawsubcubes{6.5}{2.5}{6.48}{4}
        \drawsubcubes{6.5}{2.5}{6.48}{2}
        \drawsubcubes{6.5}{4.5}{6.48}{2}
        \drawsubcubes{6.5}{2.5}{8.48}{2}
        \drawsubcubes{6.5}{4.5}{8.48}{2}
    \end{scope}
    \end{tikzpicture}}
    \begin{minipage}{0.99\textwidth}
        \centering
        \subcaption{Face-sharing interaction.}
        \label{fig:subfig3}
    \end{minipage}

    \caption{Sub-division of the source $Y$ in three dimensions.}
    \label{fig:3d combined panels}
\end{figure}
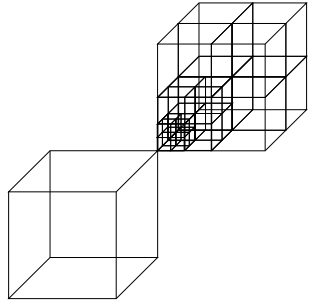
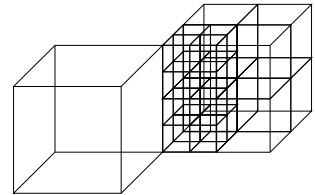
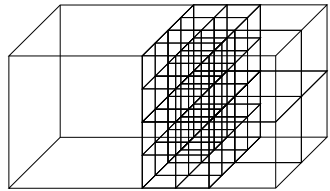

In higher dimensions, the hierarchical subdivision is carried out using an adaptive $2^d$-tree, where $Y$ is subdivided into multiple parts based on shared hyper-surfaces between the source and target domains as shown in \autoref{fig: 2d combined panels} and \autoref{fig:3d combined panels} for two- and three-dimensional interactions respectively.
\begin{itemize}
    \item At \textbf{level 1}, $Y$ is subdivided into $2^d$ equal parts $Y_{1,l}$ for $l=1:2^d$.
    \item At \textbf{level 2},  those $Y_{1,l}$'s that share a $d'$-dimensional hyper-surface with the target domain $X$ will again be subdivided into $2^d$ equal parts. This subdivision is continued up to level $\kappa$, where $\kappa = \left\lfloor \log_{2^d}n \right\rfloor $.
\end{itemize}
The hierarchical process in $d$-dimensions can be expressed as:
\[ Y = \overbrace{Y_1 \bigcup_{l=1}^{h_1} Y_{1,l}}^{\text{at level 1}} =\dotsb = \overbrace{Y_i \bigcup_{k=1}^i\bigcup_{l=1}^{h_k} Y_{k,l}}^{\text{at level }i} = \dotsb = \overbrace{Y_\kappa \bigcup_{k=1}^\kappa\bigcup_{l=1}^{h_k} Y_{k,l}}^{\text{at level }\kappa} \] where, 
 $h_k = 2^{d'k}\bkt{2^{d-d'}-1}$ for $k=1:\kappa$ and $Y_i $'s are the union those of $Y_{i,l}$'s that share a $d'$-dimensional hyper-surface with $X$ at level $i=1:\kappa$, i.e. $\displaystyle Y_i = \bigcup_{l = h_k+1}^{2^d} Y_{i,l}$ for $i = 1:\kappa$. 

\noindent We now define the random variables associated with the number of particles in each subdomain. Let:
\begin{itemize}
    \item $N_{k,l}$ ($N_{k,l}=N\bkt{Y_{k,l}}$ as defined in \autoref{equ: binomial distribution pmf}) denotes the number of particles in the sub-domain $Y_{k,l}$ for $k=1:\kappa$ and $l=1:h_k$.
    \item $M_\kappa$ ($ M_\kappa = N\bkt{Y_\kappa} $ as defined in \autoref{equ: binomial distribution pmf}) represents the number of particles in the sub-domain $Y_\kappa$, which is the final level subdivision that shares a $d'$-dimensional hyper-surface with $X$.
\end{itemize}
The total number of particles in the source domain $Y$ is the sum of the particles in $Y_\kappa$ and all the subdivisions $Y_{k,l}$'s and thus we have \[ n = M_\kappa + \sum_{k=1}^\kappa \sum_{l=1}^{h_k} N_{k,l}.\] where $n$ is the total number of particles in $Y$. Now, due to the uniform distribution, the probability of a particle being located in a specific sub-domain is proportional to the size of that sub-domain. More specifically:
\begin{itemize}
    \item The probability of a particle being located in the sub-domain $Y_{k,l}$ is given by $\mathtt{q}_{k}=\dfrac{1}{2^{dk}}$ for each $l$, where $d$ represents the dimension of the domains.
    \item The probability of a particle being in $Y_\kappa$, the sub-domain at the finest level, is $\mathtt{q}_\kappa = \dfrac{1}{2^{\bkt{d-d'}\kappa}}$, where $d'$ is the dimension of the shared hyper-surface between $X$ and $Y$.
\end{itemize}
This probabilistic structure will later allow us to estimate the rank of the random kernel matrix $K$.

\subsection{Low-Rank Matrix Construction}
The hierarchical subdivision of the source domain $Y$ allows us to generate some matrices from the kernel matrix $K$ that correspond to the interactions between the target domain $X$ and the subdivided regions of the source domain $Y$ and then efficiently approximate those matrices as a low-rank matrix by leveraging the far-field approximation provided by \autoref{lma: cheb interpolation is d dimensions}. This approach extends naturally to higher dimensions, with a one-dimensional case shown in \autoref{fig: hierarchical subdivision at level k} and \autoref{fig: hierarchical subdivision at level kappa} for pictorial clarity. The construction of the sub-matrices proceeds as follows:
\begin{itemize}
    \item The matrix $K_{k,l}$ due to the interaction between target domain $X$ and the subdivided source domain $Y_{k,l}$ for $k=1:\kappa$ and $l=1:h_k$ (the one-dimensional case is as shown in \autoref{fig: hierarchical subdivision at level k}), is given by \[\bkt{K_{{k,l}}}_{i,j} = \begin{cases}
        \mclk\bkt{x_i,y_j} &\text{where } x_i\in \text{int}(X) \;\text{ and } y_j\in \text{int}(Y_{k,l}).\\
        0 &\text{elsewhere.}
    \end{cases} \]
    \begin{figure}[H]
            \centering
            \resizebox{0.5\textwidth}{!}{
            \begin{tikzpicture}
                \draw (0,0) -- (8,0);
                \draw (0,-0.2) -- (0,0.2);
                \draw (5.1,-0.15) -- (5.1,0.15);
                \draw (4.5,-0.15) -- (4.5,0.15);
                \draw (4,-0.2) -- (4,0.2);
                \draw (8,-0.2) -- (8,0.2);
                \foreach \i in {1,...,80}{
                    \pgfmathsetmacro{\randomx}{6+2*rand}
                    \fill[blue] (\randomx-4,0) circle(0.5pt);
                }
                \foreach \i in {1,...,20}{
                    \pgfmathsetmacro{\randomx}{4.8+0.31*rand}
                    \fill[red] (\randomx,0) circle(0.45pt);
                }
                \node at (4.8,0.4) {\tiny$Y_{k,1}$};
                \node at (2,0.4) {$X$};
                \node at (0,-0.4) {$-l$};
                \node at (3.99,-0.4) {$0$};
                \node at (8,-0.4) {$l$};

                \pgfmathsetmacro{\xx}{10}
                \pgfmathsetmacro{\yy}{-2}
                
                \fill[gray!30] (\xx+1,\yy+0) rectangle (\xx+2.4,\yy+4);
                
                \draw (\xx+0,\yy+0) rectangle (\xx+4,\yy+4);
                
                \node at (\xx-0.6,0) {$K_{k,1}=$};
                \draw[white, line width=1.2pt] (\xx+0.16,\yy+0) -- (\xx+3.85,\yy+0);
                \draw[white, line width=1.2pt] (\xx+0.16,\yy+4) -- (\xx+3.85,\yy+4);    
            \end{tikzpicture}}
            \caption{ The target domain $X$ and the subdivided source domain $Y_{k,1}$ at level $k$ with the corresponding matrix $K_{k,1}$.}
            \label{fig: hierarchical subdivision at level k}
        \end{figure}
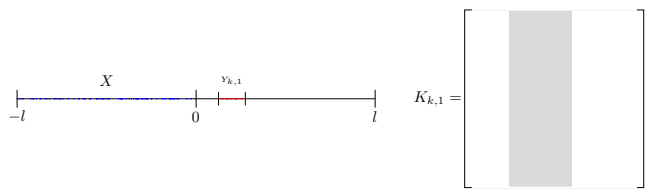
    \item The matrix $K_{\kappa}$ due to the interaction between target domain $X$ and the subdivided source domain $Y_{\kappa}$ (the one-dimensional case is as shown in \autoref{fig: hierarchical subdivision at level kappa}), is given by \[\bkt{K_\kappa}_{i,j} = \begin{cases}
        \mclk\bkt{x_i,y_j} &\text{where } x_i\in \text{int}(X) \;\text{ and } y_j\in \text{int}(Y_{\kappa}).\\
        0 &\text{elsewhere.}\end{cases} \]
    \begin{figure}[H]
        \centering
        \resizebox{0.5\textwidth}{!}{
        \begin{tikzpicture}
            \draw (0,0) -- (8,0);
            \draw (0,-0.2) -- (0,0.2);
            \draw (4.25,-0.15) -- (4.25,0.15);
            \draw (4,-0.2) -- (4,0.2);
            \draw (8,-0.2) -- (8,0.2);
            \foreach \i in {1,...,80}{
                \pgfmathsetmacro{\randomx}{6+2*rand}
                \fill[blue] (\randomx-4,0) circle(0.5pt);
            }
            \foreach \i in {1,...,10}{
                \pgfmathsetmacro{\randomx}{4.12+0.1*rand}
                \fill[red] (\randomx,0) circle(0.45pt);
            }
            \node at (4.15,0.4) {\tiny$Y_{\kappa}$};
            \node at (2,0.4) {$X$};
            \node at (0,-0.4) {$-l$};
            \node at (3.99,-0.4) {$0$};
            \node at (8,-0.4) {$l$};
            
            \pgfmathsetmacro{\xx}{10}
            \pgfmathsetmacro{\yy}{-2}
            
            \fill[gray!30] (\xx+0,\yy+0) rectangle (\xx+0.3,\yy+4);
            
            \draw (\xx+0,\yy+0) rectangle (\xx+4,\yy+4);
            
            \node at (\xx-0.55,0) {$K_\kappa=$};
            \draw[white, line width=1.2pt] (\xx+0.16,\yy+0) -- (\xx+3.85,\yy+0);
            \draw[white, line width=1.2pt] (\xx+0.16,\yy+4) -- (\xx+3.85,\yy+4);    
        \end{tikzpicture}}
        \caption{ The target domain $X$ and the subdivided source domain $Y_{\kappa}$ at level $\kappa$, with corresponding  matrix $K_{\kappa}$.}
        \label{fig: hierarchical subdivision at level kappa}
    \end{figure}
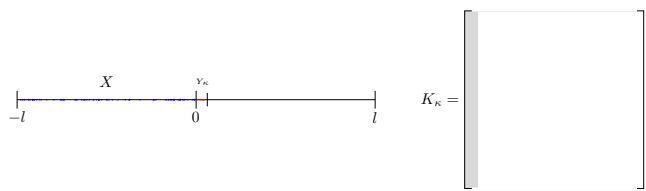
\end{itemize}
Hence, we can write the kernel matrix $K$ as the following sum \begin{equation}\label{equ: K as a sum of matrices}
    K=\sum_{k=1}^\kappa\sum_{l=1}^{h_k}K_{k,l} + K_\kappa
\end{equation}
Now, if $\Tilde{K}_{k,l}$ and $\Tilde{K}_\kappa$ be the approximations of $K_{k,l}$ and $K_\kappa$ respectively, then the approximation of the kernel matrix $K$ can be written as follows \begin{equation}\label{equ: tilde K as a sum of matrices}
    \Tilde{K}=\sum_{k=1}^\kappa\sum_{l=1}^{h_k}\Tilde{K}_{k,l} + \Tilde{K}_\kappa
\end{equation}
Here, the matrix $\Tilde{K}$ is a low-rank approximation of the kernel matrix $K$.

\subsection{Estimating Random Rank: The Random Variable \texorpdfstring{$\mclr$}{R}} \label{subsec: Estimating Random Rank: The Random Variable R}
With all the necessary components now in place, we are ready to define the random variable $\mclr$, which encapsulates the complexity of the rank behavior of the approximated matrix $\Tilde{K}$ of he kernel matrix $K$ due to the interaction between the source domain $Y$ and the target domain $X$ through the kernel function $\mclk$. Here, $\mclr$ depends on how well the matrices $K_{k,l}$ and $K_\kappa$ are approximated, which in turn depends on the geometry of the domains and the accuracy of the approximation.

The random variable $\mclr$ is formulated using the matrices $K_\kappa$ and $K_{k,l}$ for $k=1:\kappa$ and $l=1:h_k$, incorporating their ranks into the following expression:
\begin{equation}\label{equ: random rank in dD}
    \boxed{\mclr = \sum_{k=1}^\kappa\sum_{l=1}^{h_k} \mclz_n^{k,l} + M_\kappa, \quad \text{where } \mclz_n^{k,l}=\min\bkct{N_{k,l},p} } 
\end{equation}
Here’s a breakdown of the terms:
\begin{itemize}
    \item The random variable $\mclz_n^{k,l}$ serves as an upper bound of the rank of the random matrix $ \Tilde{K}_{k,l} $. Its value is influenced by the approximation accuracy $ \delta $ and depends on both the number of particles $N_{k,l}$ in $Y_{k,l}$ and the rank truncation parameter $p$, which is determined by \autoref{lma:farfield}.
    \item The random variable $M_\kappa$ represents an upper bound of the rank of the matrix $\Tilde{K}_{\kappa}$, capturing the interaction of $X$ and $Y_{\kappa}$.
\end{itemize}

\subsection{Results from Hierarchical Subdivision and Random Variable Definitions}
With the hierarchical subdivision of the domains and the definition of the corresponding random variables $N_{k,l}$, $M_\kappa$ and $\mclz_n^{k,l}$, we are now in a position to state some immediate results which provide some insight into the expected behavior, variance, and covariance of the random ranks associated with the matrices as the number of particles $n$ grows large.

\begin{myresult}
   For any fixed positive integer $k \leq \kappa$ and for all $l=1:h_k$, the limit of the expected random rank $\mclz_n^{k,l}$ of the corresponding matrix $\Tilde{K}_{k,l}$ is $p$, i.e., \[ \bbeb{\mclz_n^{k,l}} \to p,\quad  \text{as } n\to \infty.\]
\end{myresult}
\begin{proof} From \autoref{equ: expectation of Znab compact form}, we have  \begin{equation} \label{equ: EZnk} \bbeb{\mclz_n^{k,l}} = p - \sum_{i=0}^p i \,\bbpb{N_{k,l}=i} \end{equation}
Clearly, $\bbeb{\mclz_{n}^{k,l}}\leq p$. Now for $i=1:p$, $\exists\,  \Xi >0$ ($ \Xi $ independent of $n$)  such that \begin{align}\label{equ: Prob Nkl upper bound}
    \bbpb{N_{k,l} = i} &\leq  \int\limits_{a_{i,k}^{(n)}}^{b_{i,k}^{(n)}} \frac{1}{\sqrt{2\pi}}e^{-t^2/2}dt + \frac{ \Xi }{\sqrt{n}}
\end{align}
Here, $a_{i,k}^{(n)}= \frac{i-0.5 - \mu_k}{\sigma_k} $ and $b_{i,k}^{(n)} = \frac{i+0.5 - \mu_k}{\sigma_k} $, where $\mu_k=n\mathtt{q}_k$, $\sigma_k=\sqrt{n\mathtt{q}_k\bkt{1-\mathtt{q}_k}}$. Now, for some $c_{i,k}^{(n)}$ such that $\Tilde{a}_{i,k}^{(n)} \leq c_{i,k}^{(n)} \leq \Tilde{b}_{i,k}^{(n)}$, we have \begin{align*}
    \bbeb{\mclz_n^{k,l}} &\geq p - \sum_{i=0}^p i \frac{1}{\sqrt{2\pi}}e^{-\bkt {c_{i,k}^{(n)}}^2/2} + \frac{p(p+1)}{2} \frac{ \Xi }{\sqrt{n}}
\end{align*}
Now, as $n\to\infty$, $c_{i,k}^{(n)} \to - \infty$  and hence, we have $\displaystyle \lim_{n\to\infty} \bbe{\bkbt{\mclz_n^{k,l}} } = p$.
\end{proof}

\begin{myresult}
    For any fixed positive integer $k \leq \kappa$ and for all $l=1:h_k$, the limit of the variance of random rank $\mclz_n^{k,l}$ of the corresponding matrix $\Tilde{K}_{k,l}$ is $0$, i.e., \[ \var{\mclz_n^{k,l}} \to 0,\quad  \text{as } n\to \infty.\]
\end{myresult}
\begin{proof}
    Using \autoref{equ: variance Znab}, we have an upper bound the variance of $\mclz_n^{k,l}$ as
    \begin{align*}
        \var{ \mclz_n^{k,l} } &= \sum_{i=0}^{p} i^2\bbp\bkt{N_{k,l}=i} - \bkt{\sum_{i=0}^{p} i \bbp\bkt{N_{k,l}=i}}^2 \leq \sum_{i=0}^{p} i^2\bbp\bkt{N_{k,l}=i}
    \end{align*}
    As $k$ and $p$ are fixed, there exists $n_0\in\bbn$ such that for all $n>n_0$, we have  $p < \mu_{k,l} = n\mathtt{q}_{k} $ and the probability $\bbp\bkt{N_{k,l}=i}$ increases as $i$ increases up to $p$ and thus we have \begin{equation}\label{equ: variance Znkl with probability}
        \var{ \mclz_n^{k,l} }\leq \frac{p(p+1)(2p+1)}{6}\bbp\bkt{N_{k,l}=p}
    \end{equation}
    Now using \autoref{equ: Prob Nkl upper bound}, we have upper bound of $\var{ \mclz_n^{k,l} }$ as
    \begin{align}\label{equ: upper bound of varZnk}
        \var{ \mclz_n^{k,l} }&\leq \frac{p(p+1)(2p+1)}{6} \bkt{\frac{1}{\sqrt{2\pi}}e^{-\bkt{c_{p,k}^{(n)}}^2/2} + \frac{ \Xi }{\sqrt{n}} }, \quad\text{where } a_{p,k}^{(n)} \leq c_{p,k}^{(n)} \leq b_{p,k}^{(n)}. \notag\\
                    &\leq \frac{p(p+1)(2p+1)}{6} \bkt{\frac{1}{\sqrt{2\pi}}e^{-\bkt{c_{p,k}^{(n)}}^2/2} + \frac{ \Xi }{\sqrt{n}} }
    \end{align}
    Now, $a_{i,k}^{(n)} \to - \infty$ as $n\to\infty$ and hence, we have $\displaystyle \lim_{n\to\infty} \var{\mclz_n^{k,l} } = 0 $.
\end{proof}

\begin{myresult}
    For any $k_1,k_2$ and $l_1,l_2$ such that $ k_1,k_2=1:\kappa $ and $l_1,l_2=1:h_k$, the limit of the covariance of the random ranks $\mclz_n^{k_1,l_1}$ and $\mclz_n^{k_2,l_2}$ corresponding to the matrices $\Tilde{K}_{k_1,l_1}$ and $\Tilde{K}_{k_2,l_2}$ respectively is 0, i.e., \[ \cov{\mclz_n^{k_1,l_1},\mclz_n^{k_2,l_2}} \to 0,\quad  \text{as } n\to \infty. \]
\end{myresult}

\begin{proof}
    \textbf{Case I:} $k_1 \neq k_2 $ and for any $l=1:h_k$. The result follows from the inequality $\abslt{\cov{ \mclz_n^{k_1,l},\mclz_n^{k_2,l} }} \leq \sqrt{ \var{\mclz_n^{k_1,l} } \var{\mclz_n^{k_2,l} } }$ and for fixed $k_1,k_2$; $\var{\mclz_n^{k_i,l}} \to 0,\text{ as } n\to \infty$ where $i=1,2$.\\
    \textbf{Case II:} $k_1 = k_2 = k$ (say) and $l_1 \neq l_2 $. This case will arise in all $d$-dimensional settings except for $d = 1$. Now as for some fixed $k$, the sub hyper-cubes $Y_{k,l_1}$ and $Y_{k,l_2}$ of $Y$ are identical but a shift only, the distribution of $\mclz_n^{k,l_1}$ and $\mclz_n^{k,l_2}$ defined in \autoref{equ: minNkp pfm} and their joint distribution defined in \autoref{equ: minZniZnj pdf} will be the same for all $l_1,l_2=1:h_k$ and thus we have $\displaystyle \cov{\mclz_n^{k,l_1},\mclz_n^{k,l_2}} = \var{\mclz_n^{k,l}},\; \text{ for any } l=1:h_k. $
    Now as $\var{\mclz_n^{k,l}} \to 0,$ as $ n\to \infty$, we can conclude that $\displaystyle \cov{\mclz_n^{k,l_1}, \mclz_n^{k,l_2}} \to 0, \text{ as } n\to \infty. $
\end{proof}
\begin{myremark}
    The above result is also true if one of the random variables is replaced by $M_\kappa$, i.e., for some fixed $k\neq \kappa$ \[ \cov{\mclz_n^{k,l}, M_\kappa} \to 0,\quad  \text{as } n\to \infty. \]
\end{myremark}


\begin{myremark}
    As $n\to \infty$, the random variables $ \mclz_n^{k_1,l_1} $ and $ \mclz_n^{k_2,l_2} $ (for some $k_1,k_2$ and $l_1, l_2$) converge to the $p$ with no variability. Thus, in the limit, $ \mclz_n^{k,l_1} $ and $ \mclz_n^{k,l_2} $ both behave like deterministic variables that always take the value $p$.
\end{myremark}

Now, to establish a solid foundation for the proofs (in \autoref{sec: Detailed Proof of the Results}) of our main theorems, we explore a few key lemmas that will be crucial in understanding the behavior of the random variable $\mclr$.
\begin{mylma} \label{lma: sum of Znk1s is constant}
    For any large value of $n$, if $\kappa = \lfloor \log_{2^d}n\rfloor $ then for some fixed $l$ the sum of variances $\displaystyle \sum_{k=1}^\kappa \var{\mclz_n^{k,l}} $ is bounded by $C$ such that $C\in \mclo{\log\log\log n}$. 
\end{mylma}
\begin{proof}
    Without Loss of Generality, let us take $l=1$ (as \autoref{equ: upper bound of varZnk} is independent of $l$). Now, there exist $\tilde{k}$ such that the term $\bbpb{N_{k,1}=p}$ in \autoref{equ: variance Znkl with probability} is sufficiently small for all $k\leq \Tilde{k}$. More precisely, for any small $\varepsilon>0$, we have $\bbpb{N_{k,1}=p}<\varepsilon$ for all $k\leq \Tilde{k}$, where $\Tilde{k}= \left\lfloor \log_{2^d}\bkt{1 + \frac{n+1-2p}{M^2 +2p-1}} \right\rfloor $, $M \in \mclo{\log\bkt{{1}/{\varepsilon}}}$ and also $\kappa - \Tilde{k}< \Omega + \log_{2^d}\bkt{{M^2 +2p-1} }$, where $\Omega$ is some constant (independent of $n$). Now if we choose $\varepsilon = \frac{6c_1}{\kappa p (p+1)(2p+1)} $, then \begin{align*}
        \sum_{k=1}^{\Tilde{k}}\var{ \mclz_n^{k,1} } \leq \sum_{k=1}^{\Tilde{k}} \frac{p(p+1)(2p+1)}{6}\bbp\bkt{N_{k,1}=p} < c_1
    \end{align*}
    where $c_1$ is a constant independent of $n$. Now as $\bkt{\kappa -\Tilde{k}} \in \mclo{\log\log\log n}$ and $\var{\mclz_n^{k,1}}$ is bounded by $p(p+1)(2p+1)/6$ (from \autoref{equ: variance Znkl with probability}), we can conclude that there exist $C$ such that \[\sum_{k=1}^\kappa \var{\mclz_n^{k,1}}\leq C, \quad \text{where $C\in \mclo{\log\log\log n}$.} \] 
\end{proof}

\begin{mylma}\label{lma: sum of ZniZnjs is constant}
    For any large value of $n$, if $\kappa = \lfloor \log_{2^d}n\rfloor $ then for some fixed $l$ the sum of covariances \[\sum_{k_1=1}^\kappa\sum_{k_2=k_1+1}^\kappa\cov{\mclz_n^{k_1,l},\mclz_n^{k_2,l}}\leq C,\quad \text{where $C\in \mclo{(\log\log\log n)^2} $.}\] 
\end{mylma}
\begin{proof}
    Without Loss of Generality, let us take $l=1$ as for some fixed $k_1,k_2$, then the upper bound of $\cov{\mclz_n^{k_1,l},\mclz_n^{k_2,l}}$ defined in \autoref{equ: cov(ZniZnj)} can be obtained as \begin{align*}
        \cov{\mclz_n^{k_1,l}, \mclz_n^{k_2,l}} \leq \sum_{r=0}^p\sum_{m=0}^p rs\bbpb{N_{k_1,1}=p_1,N_{k_2,1}=p_2},
    \end{align*}
    for some $p_1,p_2$ such that  $1\leq p_1,p_2 \leq p$. Now for any small enough $\varepsilon>0$, we have \[\bbpb{N_{k_1,1}=p_i,N_{k_2,1}=p_j}<\varepsilon,\; \text{ for all }k_1,k_2\leq \Tilde{k}.\] Here, $\Tilde{k}$ will take a similar form like in \autoref{lma: sum of Znk1s is constant} and hence in a similar approach, we can conclude that there exists $C$ such that \[ \sum_{k_1=1}^\kappa\sum_{k_2=k_1+1}^\kappa\cov{\mclz_n^{k_1,l},\mclz_n^{k_2,l}}\leq C, \quad \text{where $C\in \mclo{(\log\log\log n)^2}$}.\]
\end{proof}
\noindent\textbf{Note: } A detailed version of the proofs of the \autoref{lma: sum of Znk1s is constant} and \autoref{lma: sum of ZniZnjs is constant} can be found in \autoref{app: proof of lemma: sum of Zni1s is constant and lemma: sum of ZniZnjs is constant}.

\begin{mylma}\label{lma: sum of ZniMkappa is constant}
    For any value of $n$, if $\kappa = \lfloor \log_{2^d}n\rfloor $ then for some fixed $l$ the sum of covariances \[ \sum_{k=1}^\kappa\cov{\mclz_n^{k,l}, M_\kappa} \in \mclo{n^{d'/d}\log\log\log n}. \]
\end{mylma}
\begin{proof}
    An upper bound of $ \cov{\mclz_n^{k,l}, M_\kappa} $ (defined in \autoref{equ: cov(Znk_Mkappa)}) can be obtained as \[ \cov{ \mclz_n^{k,l} , M_\kappa} \leq p^2\sum_{m=0}^{n} m \bbpb{N_{k,1}=p_1,M_\kappa=m} \] for some $0\leq p_1 \leq p$.  Using conditional expectation, the  above expression can be rewritten as \[ \cov{ \mclz_n^{k,l}, M_\kappa} \leq p^2 \,\bbeb{M_\kappa|N_{k,1}=p_1}\bbpb{N_{k,1}=p_1} \]
    Now as $\bkt{M_\kappa|N_{k,1}=p_1} \sim \textbf{Binomial}\bkt{n-p_1,\dfrac{\mathtt{q}_{\kappa}}{1-\mathtt{q}_k}}$, we have $\bbeb{M_\kappa|N_k=p_1} = \bkt{n-p_1}\dfrac{\mathtt{q}_\kappa}{1-\mathtt{q}_k}$. Hence, \[ \cov{ \mclz_n^{k,l} , M_\kappa} \leq p^2 \frac{n-p_1}{2^{\bkt{d-d'}\kappa}} \frac{2^{dk}}{2^{dk}-1} \bbpb{N_{k,1}=p_1} \]
    Now as $\dfrac{n-p_1}{2^{d\kappa}}\in \mclo{1}$ and $\dfrac{2^{dk}}{2^{dk}-1}$ are bounded by 2, there exist $c$ such that \[ \cov{ \mclz_n^{k,l} , M_\kappa} \leq c\,n^{d'/d}\, \bbpb{N_{k,1}=p_1}.\]
    Now similarly like \autoref{lma: sum of Znk1s is constant}, we can prove that \[ \sum_{k=1}^\kappa \cov{\mclz_n^{k,l}, M_\kappa} \in \mclo{n^{d'/d}\log\log\log n} \]
\end{proof}
\begin{myremark}\label{rmk: covariance is bounded by constant}
    For the vertex-sharing case, i.e., for $d'=0$ in any dimension, there exists $C$ such that \[ \sum_{k=1}^\kappa \cov{\mclz_n^{k,l}, M_\kappa} \leq  C, \quad \text{where $C\in \mclo{\log\log\log n}$}. \]
\end{myremark}

\section{Proof of Theorems}\label{sec: Detailed Proof of the Results}
\begin{figure}[H]
    \centering
    \subfloat[\footnotesize one-dimensional case]{%
        \includegraphics[width=0.45\textwidth]{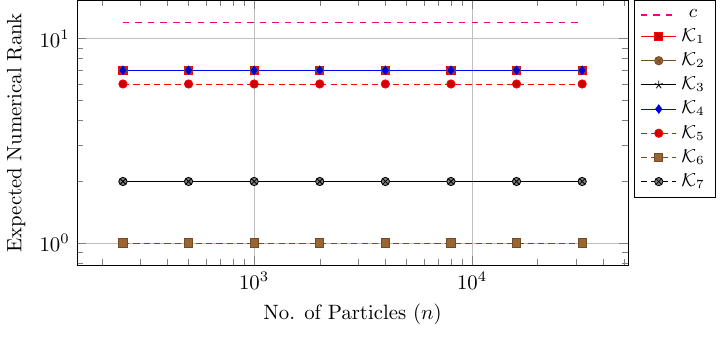}
    }
    \subfloat[\footnotesize two-dimensional case]{%
        \includegraphics[width=0.45\textwidth]{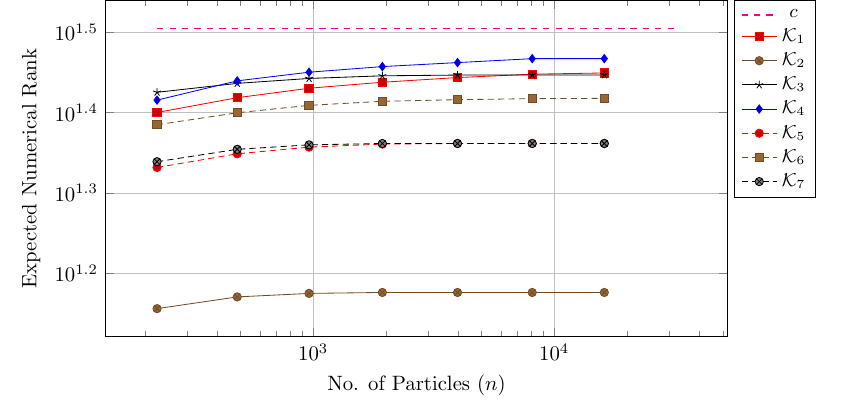}
    }\\ \vskip 0.2cm
    \subfloat[\footnotesize three-dimensional case]{%
        \includegraphics[width=0.45\textwidth]{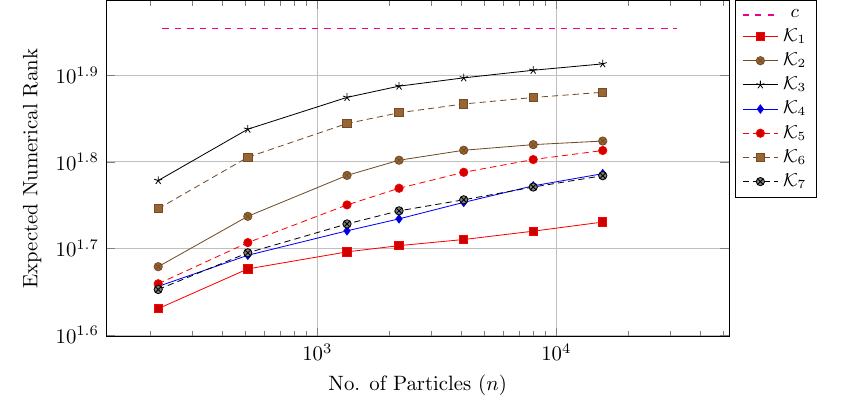}
    }
    \caption{Expected Numerical Rank Growth of different Kernel matrices for one-, two- and three-dimensional far-field interactions.}
    \label{plot: far-field rank growth}
\end{figure}

In this section, we present the proofs of the theorems stated in \autoref{sec: Main Results}. In these proofs, we will use the well-known result that for any $d$-dimensional case, the rank kernel matrix, due to the far-field interaction [\!\!\cite{KHAN2024112786} Lemma 5.1] with arbitrarily distributed particles, is bounded by a constant, which dependents only on the desired accuracy $ \delta $, as validated by the results\footnote{The experiments were repeated 2000 times for matrices of size less than 8100 and 500 times for matrices of size greater than 8100 in order to obtain the sample data of Numerical Ranks $\bkt{\text{with } \delta =10^{-12}}$ in one-, two- and three-dimensions with all possible interactions of domains.} shown in \autoref{plot: far-field rank growth}.

\subsection{Expected Rank Growth}
In this subsection, we prove the expected growth of the random rank $\mclr$ for any $d'$-dimensional hyper-surface sharing domains in $d$-dimensions.


\begin{customproof}{\autoref{thm: expected rank growth in dD}}  
    The expected value of $\mclr$ (defined in \autoref{equ: random rank in dD}) is given by 
    \begin{align*}
        \bbeb{\mclr} &= \sum_{k=1}^\kappa\sum_{l=1}^{h_k} \bbeb{\mclz_n^{k,l}} + \bbeb{M_\kappa}
    \end{align*}
    Now from \autoref{equ: EZnk}, the expectation of $\mclz_n^{k,l}$ is given by 
    \[ \bbeb{\mclz_n^{k,l}} = p - \sum_{i=0}^p i\,\bbpb{N_{k,l} = i} \]
    Thus we have \begin{align*}
        \bbeb{\mclr} &= \sum_{k=1}^\kappa h_k \bkt{p - \sum_{i=0}^p i\,\bbp\bkt{N_{k,l} = i}} + \bbeb{M_\kappa}\\
                         &= \sum_{k=1}^\kappa h_k p - \sum_{k=1}^\kappa\sum_{i=0}^p h_ki\,\bbp\bkt{N_{k,l} = i} + n^{d'/d} \frac{n}{2^{d\kappa}}
    \end{align*} 
    We now encounter the cases that $d'$ being zero and non-zero separately below.
    \begin{enumerate}[(i)]
        \item \underline{\bf $\mathbf{d'=0}$}: In this case the value of $h_k$ remains constant and given by $h_k=\bkt{2^d-1}$ for all $k = 1 : \kappa$ and hence the expectation is given by 
        \begin{align*}
            \bbeb{\mclr} &=  \bkt{2^d-1}\sum_{k=1}^\kappa p - \bkt{2^d-1}\sum_{k=1}^\kappa\sum_{i=0}^p i\,\bbp\bkt{N_{k,l} = i} + \frac{n}{2^{d\kappa}}\\
                        &\leq \bkt{2^d-1}\kappa p + \frac{n}{2^{d\kappa}}
        \end{align*}    
        \noindent Now, $\kappa = \left\lfloor \log_{2^d} n \right\rfloor $ and $\dfrac{n}{2^{d\kappa}} \in \mclo{1}$. Hence, for the fixed dimension $d$, we can conclude that 
        \[\bbe\bkbt{\mclr}\in \mathcal{O}\bkt{p\log_{2^d}n}.\]
        \item  \underline{\bf $\mathbf{d' \neq 0}$}: In this case, the value of $h_k$ does not remain constant, but with some simple algebra, we have the following expected value \begin{align*}
        \bbeb{\mclr} &= n^{d'/d} \frac{n}{2^{d\kappa}} + \frac{2^{d'}}{2^{d'}-1}\bkt{2^{d-d'}-1}\bkt{2^{d'\kappa} - 1}p - \sum_{k=1}^\kappa\sum_{i=0}^p h_ki\,\bbp\bkt{N_{k,l} = i}\\
                    &\leq n^{d'/d} \frac{n}{2^{d\kappa}} + \frac{2^d - 2^{d'}}{2^{d'}-1}\bkt{n^{d'/d} - 1}p 
    \end{align*}
    Hence, we can conclude that if $d$-dimensional source and targets shares a $d'$-dimensional hyper-surface, then  $\bbeb{\mclr} \in \mathcal{O}\bkt{pn^{d'/d}} $. 
    \end{enumerate}
\end{customproof}

\begin{figure}[H]
    \centering
    \subfloat[\footnotesize one-dimensional case]{%
        \includegraphics[width=0.45\textwidth]{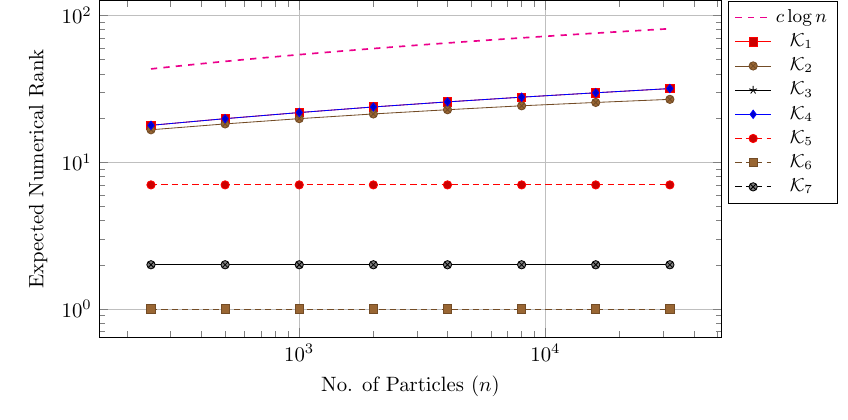}
    }
    \subfloat[\footnotesize two-dimensional case]{%
        \includegraphics[width=0.45\textwidth]{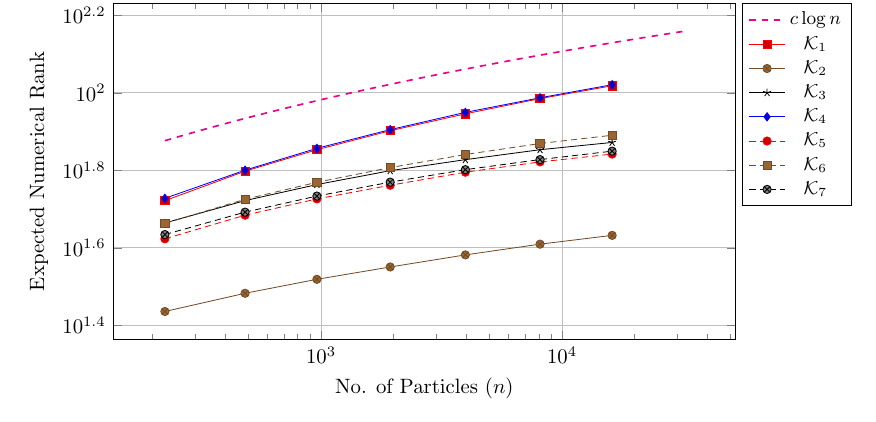}
    }\\ \vskip 0.2cm
    \subfloat[\footnotesize three-dimensional case]{%
        \includegraphics[width=0.45\textwidth]{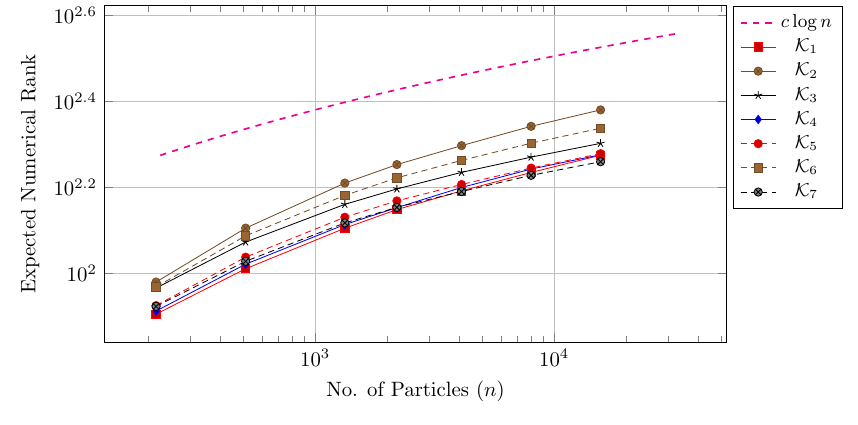}
    }
    \caption{Expected Numerical Rank Growth of Kernel matrix for different kernels for one-, two-, and three-dimensional vertex sharing interaction.}
    \label{plot: vertex-sharing rank growth}
\end{figure}

\begin{figure}[H]
    \centering
    \subfloat[\resizebox{0.35\textwidth}{!}{two-dimensional edge-sharing interaction}]{%
        \includegraphics[width=0.45\textwidth]{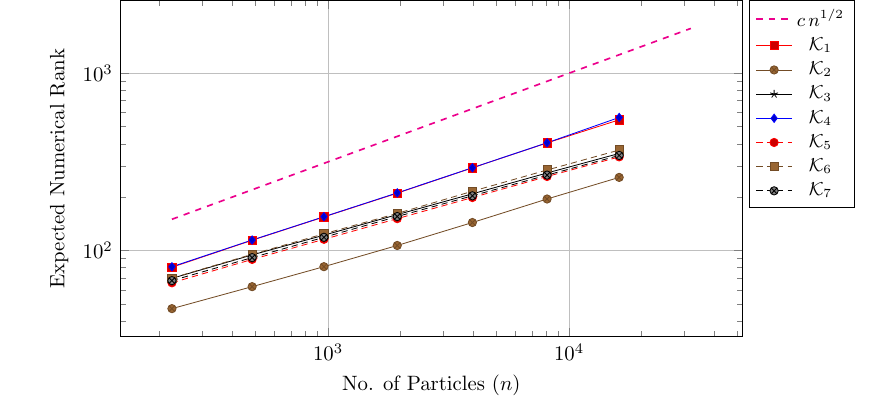}
    }\\ \vskip 0.2cm
    \subfloat[\resizebox{0.35\textwidth}{!}{three-dimensional edge-sharing interaction}]{%
        \includegraphics[width=0.45\textwidth]{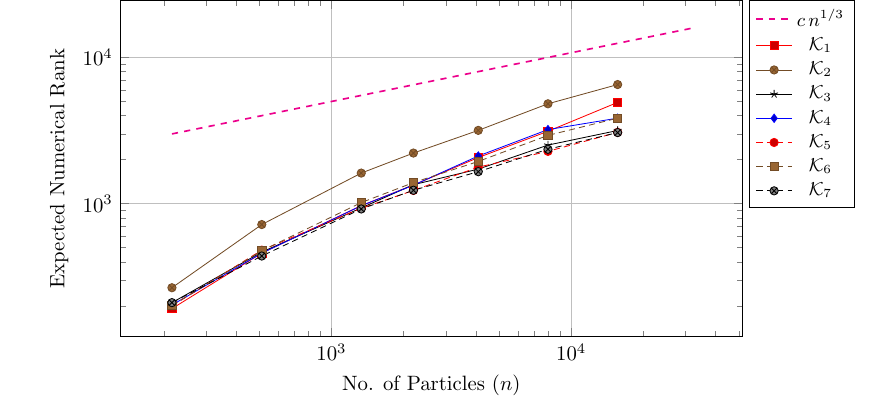}
    }
    \subfloat[\resizebox{0.35\textwidth}{!}{three-dimensional face-sharing interaction}]{%
        \includegraphics[width=0.45\textwidth]{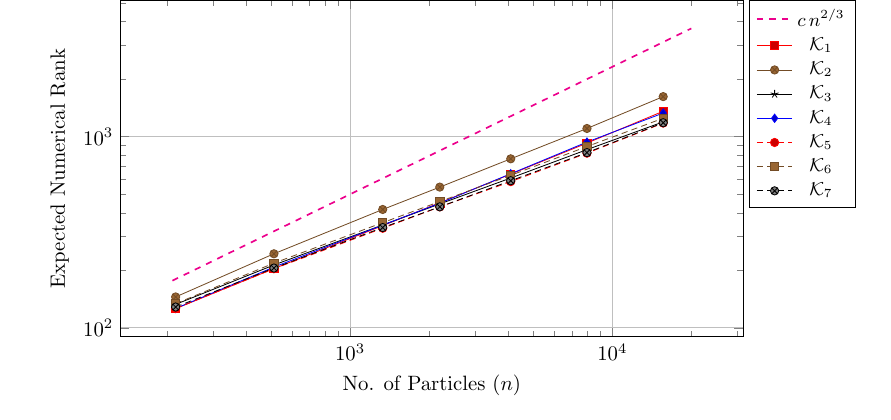}
    }
    \caption{Expected Numerical Rank Growth of Kernel matrix for different kernels for two- and three-dimensional edge sharing and three-dimensional face sharing interaction.}
    \label{plot: edge-sharing rank growth}
\end{figure}


\subsection{Growth in Variance of \texorpdfstring{$\mclr$}{R}}
In this subsection, we will prove the \autoref{thm: variance of rank in dD} on the variance growth of the random rank $\mclr$ incorporating all possible interactions in $d$-dimensions as stated in \autoref{sec: Main Results}. 

\begin{customproof}{\autoref{thm: variance of rank in dD}}
   The variance of the random variable $\mclr$ (defined in \autoref{equ: random rank in dD}) is given by \begin{equation} \label{equ: variance in dD}
        \var{\mclr} =  \sum_{k=1}^\kappa\sum_{l=1}^{h_k} \var{\mclz_n^{k,l}} + \var{M_\kappa} + 2 \sum_{k=1}^\kappa\sum_{l=1}^{h_k} \cov{\mclz_n^{k,l}, M_\kappa} + 2 \sum\limits_{\substack{\text{either } k_1< k_2\\ \text{or if } k_1 = k_2 \\ \text{then } l_1< l_2}} \cov{\mclz_n^{k_1,l_1}, \mclz_n^{k_2,l_2}}       
   \end{equation}
Now, after some simple algebraic manipulation, \autoref{equ: variance in dD} can be written as 
\begin{align}
    \var{\mclr} 
       &= \var{M_\kappa} + 2 \sum_{k=1}^\kappa h_k \cov{\mclz_n^{k,l}, M_\kappa}  + 2 \sum_{k_1=1}^\kappa\sum_{k_2=k_1+1}^\kappa h_{k_1}h_{k_2} \cov{\mclz_n^{k_1,l}, \mclz_n^{k_2,l}}  +  \sum_{k=1}^\kappa h_{k}\bkt{2h_k-2l+1} \var{\mclz_n^{k,l}} \notag\\
       &\leq \var{M_\kappa} + 2 h_\kappa \sum_{k=1}^\kappa \cov{\mclz_n^{k,l}, M_\kappa}  + 2 h_\kappa^2 \sum_{k_1=1}^\kappa\sum_{k_2=k_1+1}^\kappa \cov{\mclz_n^{k_1,l}, \mclz_n^{k_2,l}}  + 3h_\kappa^2 \sum_{k=1}^\kappa \var{\mclz_n^{k,l}} \notag\\
       &\leq 2^{d'\kappa}\frac{n}{2^{d\kappa}} +  2^{d'\kappa+1}\bkt{2^{d-d'}-1} \sum_{k=1}^\kappa \cov{\mclz_n^{k,l}, M_\kappa}  + 2^{2d'\kappa+1}\bkt{2^{d-d'}-1}^2 \sum_{k_1=1}^\kappa\sum_{k_2=k_1+1}^\kappa \cov{\mclz_n^{k_1,l}, \mclz_n^{k_2,l}} \notag\\ &\hskip 2.3cm +  3\times 2^{2d'\kappa}\bkt{2^{d-d'}-1}^2 \sum_{k=1}^\kappa \var{\mclz_n^{k,l}} \label{equ: varR upper bound expression}
\end{align}
Now we encounter the cases where $d'$ is zero and non-zero separately below.
\begin{enumerate}[(i)]
    \item \underline{\bf $\mathbf{d'=0}$}: In this case $\var{\mclr}$ is bounded by a constant which is of $ \mclo{(\log\log\log n)^2} $. i.e., $\var{\mclr}\in \mclo{(\log\log\log n)^2}$. This follows immediately from \autoref{lma: sum of Znk1s is constant}, \autoref{lma: sum of ZniZnjs is constant}, \autoref{rmk: covariance is bounded by constant}, and along with the fact $\dfrac{n}{2^{d\kappa}}\in \mclo{1}$.
    \item \underline{\bf $\mathbf{d'\neq 0}$}: Due to the increased dimensionality of hyper-surface sharing between the sources and targets, $\var{\mclr}\in \mclo{\bkt{n^{d'/d}\log\log\log n}^2}$ which follows from the \autoref{lma: sum of Znk1s is constant}, \autoref{lma: sum of ZniZnjs is constant}, \autoref{lma: sum of ZniMkappa is constant} and the fact $2^{d'\kappa}\dfrac{n}{2^{d\kappa}}\in \mclo{n^{d'/d}}$.
\end{enumerate}
\end{customproof}
\begin{myremark}
    In any dimension, for the hyper-cubes $X$ and $Y$ share a vertex, i.e., for $d'=0$ in $d$-dimensions, $\bbe\bkbt{\mclr}\in \mathcal{O}\bkt{p\log_{2^d}n}$ and $\var{\mclr}\in\mclo{(\log\log\log n)^2}$ suggests that the mean may increase (logarithmically), still the spread of the rank around its mean remains stable and does not widen much as the size of matrix increases. 
\end{myremark}

\begin{figure}[H]
    \centering
    \subfloat[\footnotesize one-dimensional case]{%
        \includegraphics[width=0.45\textwidth]{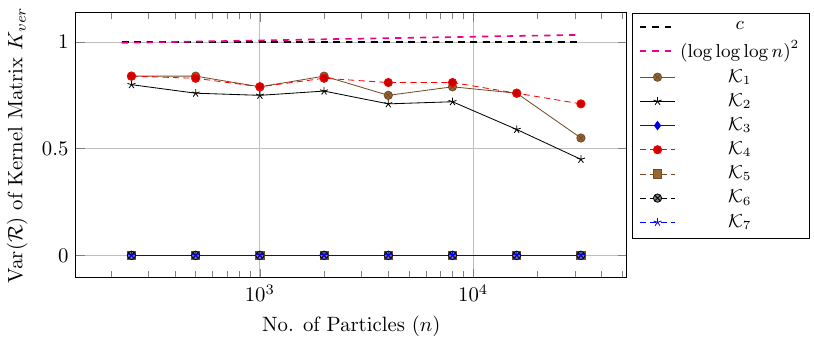}
    }
    \subfloat[\footnotesize two-dimensional case]{%
        \includegraphics[width=0.45\textwidth]{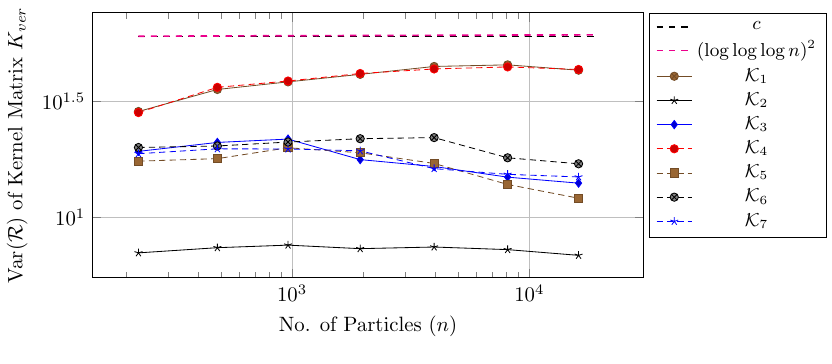}
    }\\ \vskip 0.2cm
    \subfloat[\footnotesize three-dimensional case]{%
        \includegraphics[width=0.45\textwidth]{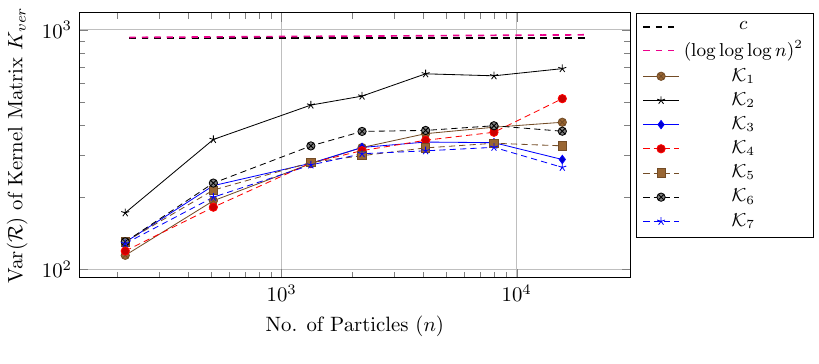}
    }
    \caption{Growth of Variance of Numerical Rank of Kernel matrix for different kernels for two- and three-dimensional vertex sharing interaction.}
    \label{plot: vertex-sharing growth of variance of rank}
\end{figure}

\begin{figure}[H]
    \centering
    \subfloat[\resizebox{0.3\textwidth}{!}{two-dimensional edge-sharing interaction}]{%
        \includegraphics[width=0.45\textwidth]{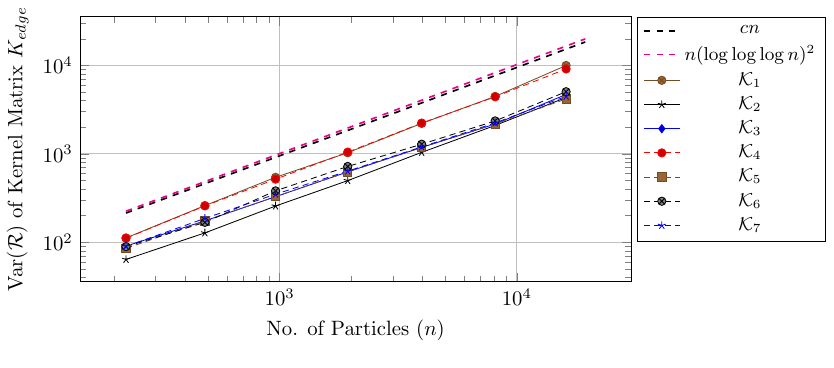}
    }\\ \vskip 0.2cm
    \subfloat[\resizebox{0.3\textwidth}{!}{three-dimensional edge-sharing interaction}]{%
        \includegraphics[width=0.45\textwidth]{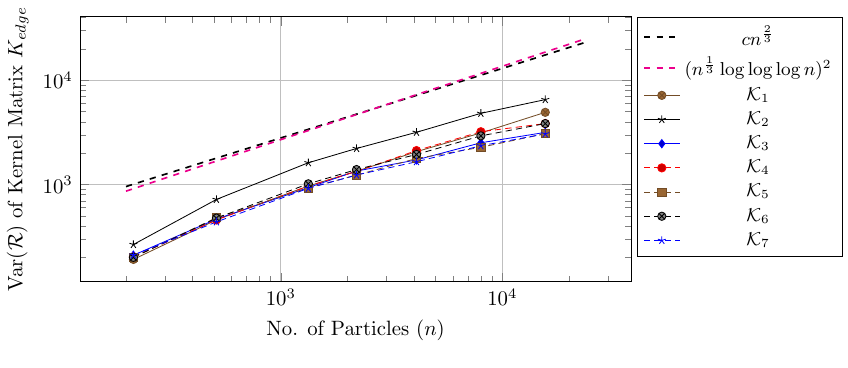}
    }
    \subfloat[\resizebox{0.3\textwidth}{!}{three-dimensional face-sharing interaction}]{%
        \includegraphics[width=0.45\textwidth]{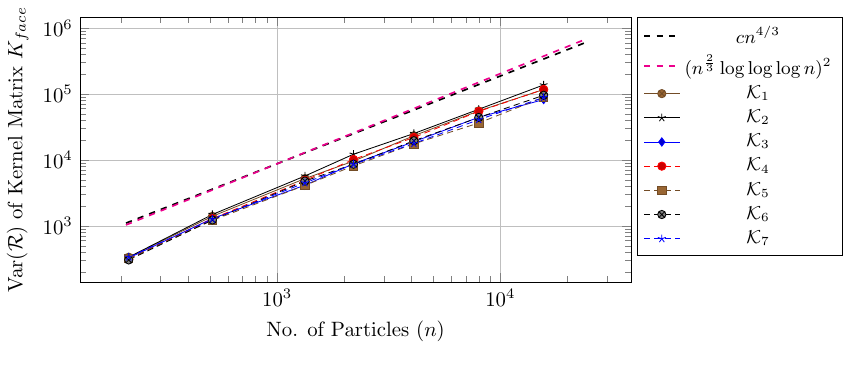}
        \label{subfig: face in 3d growth}
    }
    \caption{Growth of Variance of Numerical Rank of Kernel matrix for different kernels for two- and three-dimensional edge sharing and three-dimensional face sharing interactions.}
    \label{plot: edge-sharing growth of variance of rank}
\end{figure}

\section{Numerical Results} \label{sec: Numerical Resulsts}
In this section, we present the numerical results based on the theoretical framework outlined in \autoref{sec: Fundamental Framework and Problem Setup}. Specifically, we will explore how the expectation and variance of the random variable $\mclr$ changes depending on how the interactions between source and target domains change in one-, two-, and three-dimensions. \textit{To obtain the means and variances, we repeated the same experiment 2,000 times to collect the samples of Numerical Ranks $\bkt{\text{with }  \delta =10^{-12}}$ for matrices of size $\leq 8100$ and $500$ times for matrices of size larger than that. The experiment was conducted in \textsf{MATLAB} using parallel processing to enhance efficiency.}
\subsubsection*{Results for One-Dimensional Cases:}

\begin{table}[H]
    \centering
    \begin{minipage}{0.49\textwidth}
        \centering
        \resizebox{\textwidth}{!}{ 
        \begin{tabular}{|c| c c c c c c c c|}
            \hline
            \multirow{2}{*}{ker fun} & \multicolumn{8}{c|}{$\bbeb{\mclr}$ of Kernel Matrix $K_{far}$}  \\
            \cline{2-9}
            & $n=250$ & $n=500$ & $n=1000$ & $n=2000$  & $n=4000$& $n=8000$& $n=16000$& $n=32000$\\
            \hline
            $\mclk_1$ &7 &7 &7 &7 &7 &7 &7 &7\\ 
            $\mclk_2$ &7 &7 &7 &7 &7 &7 &7 &7\\
            $\mclk_3$ &2 &2 &2 &2 &2 &2 &2 &2\\
            $\mclk_4$ &7 &7 &7 &7 &7 &7 &7 &7\\
            $\mclk_5$ &6 &6 &6 &6 &6 &6 &6 &6 \\
            $\mclk_6$ &1 &1 &1 &1 &1 &1 &1 &1\\
            $\mclk_7$ &2 &2 &2 &2 &2 &2 &2 &2\\
            \hline
        \end{tabular}}
        \subcaption{Mean of Numerical Ranks of Data Samples}
        \label{tab: expected ranks far-field domain in 1D}
    \end{minipage}%
    \hfill 
    \begin{minipage}{0.49\textwidth}
        \centering
        \resizebox{\textwidth}{!}{ 
        \begin{tabular}{|c| c c c c c c c c|}
            \hline
            \multirow{2}{*}{ker fun} & \multicolumn{8}{c|}{$\var{\mclr}$ of Kernel Matrix $K_{far}$}  \\
            \cline{2-9}
                & $n=250$ & $n=500$ & $n=1000$ & $n=2000$  & $n=4000$& $n=8000$& $n=16000$& $n=32000$\\
            \hline
            $\mclk_1$ &0 &0 &0 &0 &0 &0 &0 &0\\ 
            $\mclk_2$ &0 &0 &0 &0 &0 &0 &0 &0\\
            $\mclk_3$ &0 &0 &0 &0 &0 &0 &0 &0\\
            $\mclk_4$ &0 &0 &0 &0 &0 &0 &0 &0\\
            $\mclk_5$ &0 &0 &0 &0 &0 &0 &0 &0\\
            $\mclk_6$ &0 &0 &0 &0 &0 &0 &0 &0\\
            $\mclk_7$ &0 &0 &0 &0 &0 &0 &0 &0\\
            \hline
        \end{tabular}}
        \subcaption{Variance of Numerical Ranks of Data Samples}
        \label{tab: variance ranks far-field domain in 1D}
    \end{minipage}
    \caption{ Random rank statistics for far-field domains in 1D.}
\end{table}

\begin{table}[H]
    \centering
    \begin{minipage}{0.49\textwidth}
        \centering
        \resizebox{\textwidth}{!}{ 
        \begin{tabular}{|c| c c c c c c c c|}
            \hline
            \multirow{2}{*}{ker fun} & \multicolumn{8}{c|}{$\bbeb{\mclr}$ of Kernel Matrix $K_{ver}$}  \\
            \cline{2-9}
            & $n=250$ & $n=500$ & $n=1000$ & $n=2000$  & $n=4000$& $n=8000$& $n=16000$& $n=32000$\\
            \hline 
            $\mclk_1$ &17.82  &19.80  &21.77 &23.77  &25.75 &27.72 &29.68  &31.70\\
            $\mclk_2$ &16.62  &18.22  &19.79 &21.30  &22.76 &24.21 &25.51  &26.78\\
            $\mclk_3$ &2.00  &2.00  &2.00 &2.00 &2.00 &2.00  &2.00  &2.00\\
            $\mclk_4$ &17.83  &19.78  &21.77 &23.77 &25.73 &27.71  &29.69  &31.60\\
            $\mclk_5$ &7.00  &7.00  &7.00 &7.00 &7.00 &7.00 &7.00  &7.00 \\
            $\mclk_6$ &1.00  &1.00 &1.00 &1.00 &1.00 &1.00 &1.00  &1.00\\
            $\mclk_7$ &2.00  &2.00  &2.00 &2.00 &2.00 &2.00 &2.00  &2.00\\
            \hline
        \end{tabular}}
        \subcaption{Mean of Numerical Ranks of Data Samples}
        \label{tab: expected ranks vertex sharing domain in 1D}
    \end{minipage}%
    \hfill 
    \begin{minipage}{0.49\textwidth}
        \centering
        \resizebox{\textwidth}{!}{ 
        \begin{tabular}{|c| c c c c c c c c|}
            \hline
            \multirow{2}{*}{ker fun} & \multicolumn{8}{c|}{$\var{\mclr}$ of Kernel Matrix $K_{ver}$}  \\
            \cline{2-9}
                & $n=250$ & $n=500$ & $n=1000$ & $n=2000$  & $n=4000$& $n=8000$& $n=16000$& $n=32000$\\
            \hline
            $\mclk_1$ &0.84 &0.84 &0.79 &0.84 &0.75 &0.79 &0.76 &0.55\\ 
            $\mclk_2$ &0.80 &0.76 &0.75 &0.77 &0.71 &0.72 &0.59 &0.45\\
            $\mclk_3$ &0 &0 &0 &0 &0 &0 &0 &0\\
            $\mclk_4$ &0.84 &0.83 &0.79 &0.83 &0.81 &0.81 &0.76 &0.71\\
            $\mclk_5$ &0 &0 &0 &0 &0 &0 &0 &0\\
            $\mclk_6$ &0 &0 &0 &0 &0 &0 &0 &0\\
            $\mclk_7$ &0 &0 &0 &0 &0 &0 &0 &0\\
            \hline
        \end{tabular}}
        \subcaption{Variance of Numerical Ranks of Data Samples}
        \label{tab: variance ranks ver sharing domain in 1D}
    \end{minipage}
    \caption{ Random rank statistics for vertex-sharing domains in 1D.}
\end{table}

\subsubsection*{Results for Two-Dimensional Cases:}

\begin{table}[H]
    \centering
    \begin{minipage}{0.49\textwidth}
        \centering
        \resizebox{\textwidth}{!}{ 
        \begin{tabular}{|c| c c c c c c c |}
            \hline
            \multirow{2}{*}{ker fun} & \multicolumn{7}{c|}{$\bbeb{\mclr}$ of Kernel Matrix $K_{far}$}  \\
            \cline{2-8}
                & $n=225$ & $n=484$ & $n=961$ & $n=1936$  & $n=3969$& $n=8100$& $n=16129$\\
            \hline
            $\mclk_1$ &25.14 &26.23 &26.95 &27.42 &27.79 &28.06 &28.16 \\ 
            $\mclk_2$ &14.32 &14.81 &14.96 &15.00 &15.00 &15.00 &15.00 \\
            $\mclk_3$ &26.64 &27.33 &27.72 &27.93 &27.99 &28.00 &28.00 \\
            $\mclk_4$ &26.05 &27.53 &28.22 &28.68 &29.01 &29.26 &29.34 \\
            $\mclk_5$ &21.47 &22.33 &22.76 &22.96 &23.00 &23.00 &23.00  \\
            $\mclk_6$ &24.27 &25.12 &25.66 &25.96 &26.08 &26.16 &26.18 \\
            $\mclk_7$ &21.83 &22.62 &22.91 &23.00 &23.00 &23.00 &23.00 \\
            \hline
        \end{tabular}}
        \subcaption{Mean of Numerical Ranks of Data Samples}
        \label{tab: expected ranks far-field domain in 2D}
    \end{minipage}%
    \hfill 
    \begin{minipage}{0.49\textwidth}
        \centering
        \resizebox{\textwidth}{!}{ 
        \begin{tabular}{|c| c c c c c c c|}
            \hline
            \multirow{2}{*}{ker fun} & \multicolumn{7}{c|}{$\var{\mclr}$ of Kernel Matrix $K_{far}$}  \\
            \cline{2-8}
                & $n=225$ & $n=484$ & $n=961$ & $n=1936$  & $n=3969$& $n=8100$& $n=16129$\\
            \hline
            $\mclk_1$ &1.32 &0.86 &0.69 &0.48 &0.40 &0.32 &0.21\\ 
            $\mclk_2$ &0.72 &0.25 &0.05 &0.00 &0.00 &0.00 &0.00\\
            $\mclk_3$ &0.82 &0.32 &0.21 &0.06 &0.01 &0.00 &0.00   \\
            $\mclk_4$ &1.83 &1.44 &0.64 &0.39 &0.31 &0.26 &0.24\\
            $\mclk_5$ &0.89 &0.53 &0.23 &0.04 &0.00 &0.00 &0.00 \\
            $\mclk_6$ &0.69 &0.71 &0.36 &0.12 &0.09 &0.13 &0.14\\
            $\mclk_7$ &0.89 &0.37 &0.09 &0.00 &0.00 &0.00 &0.00\\
            \hline
        \end{tabular}}
        \subcaption{Variance of Numerical Ranks of Data Samples}
        \label{tab: variance ranks far-field domain in 2D}
    \end{minipage}
    \caption{ Random rank statistics for far-field domains in 2D.}
\end{table}

\begin{table}[H]
    \centering
    \begin{minipage}{0.49\textwidth}
        \centering
        \resizebox{\textwidth}{!}{ 
        \begin{tabular}{|c| c c c c c c c|}
            \hline
            \multirow{2}{*}{ker fun} & \multicolumn{7}{c|}{$\bbeb{\mclr}$ of Kernel Matrix $K_{ver}$}  \\
            \cline{2-8}
            & $n=225$ & $n=484$ & $n=961$ & $n=1936$  & $n=3969$& $n=8100$& $n=16129$\\
            \hline
            $\mclk_1$ &52.70 &62.79 &71.33 &79.85 &88.28 &96.64 &104.29\\ 
            $\mclk_2$ &27.28 &30.40 &33.03 &35.54 &38.19 &40.70 &42.88\\
            $\mclk_3$ &46.16 &52.72 &57.99 &62.97 &67.24 &71.33 &74.48\\
            $\mclk_4$ &53.43 &63.17 &71.91 &80.37 &89.04 &97.11 &104.96\\
            $\mclk_5$ &42.01 &48.34 &53.25 &57.75 &62.41 &66.35 &69.50 \\
            $\mclk_6$ &46.04 &53.14 &58.78 &64.15 &69.33 &74.02 &77.74\\
            $\mclk_7$ &43.09 &49.23 &54.16 &58.83 &63.36 &67.28 &70.73\\
            \hline
        \end{tabular}}
        \subcaption{Mean of Numerical Ranks of Data Samples}
        \label{tab: expected ranks vertex sharing domain in 2D}
    \end{minipage}%
    \hfill 
    \begin{minipage}{0.49\textwidth}
        \centering
        \resizebox{\textwidth}{!}{ 
        \begin{tabular}{|c| c c c c c c c |}
            \hline
            \multirow{2}{*}{ker fun} & \multicolumn{7}{c|}{$\var{\mclr}$ of Kernel Matrix $K_{ver}$}  \\
            \cline{2-8}
                & $n=225$ & $n=484$ & $n=961$ & $n=1936$  & $n=3969$& $n=8100$& $n=16129$\\
            \hline
            $\mclk_1$ &28.54 &35.52 &38.30 &41.26 &44.62 &45.34 &43.00\\ 
            $\mclk_2$ &7.02 &7.39 &7.58 &7.32 &7.44 &7.25 &6.85\\
            $\mclk_3$ &19.22 &21.00 &21.71 &17.71 &16.54 &14.88 &14.02\\
            $\mclk_4$ &28.24 &36.28 &38.63 &41.61 &43.53 &44.42 &43.31 \\
            $\mclk_5$ &17.44 &17.88 &19.92 &18.96 &17.05 &13.84 &12.07\\
            $\mclk_6$ &19.97 &20.27 &21.08 &21.79 &22.03 &18.03 &16.99\\
            $\mclk_7$ &18.78 &19.71 &19.68 &19.32 &16.19 &15.30 &14.91\\
            \hline
        \end{tabular}}
        \subcaption{Variance of Numerical Ranks of Data Samples}
        \label{tab: variance ranks ver sharing domain in 2D}
    \end{minipage}
    \caption{ Random rank statistics for vertex-sharing domains in 2D.}
\end{table}

\begin{table}[H]
    \centering
    \begin{minipage}{0.49\textwidth}
        \centering
        \resizebox{\textwidth}{!}{ 
        \begin{tabular}{|c| c c c c c c c|}
            \hline
            \multirow{2}{*}{ker fun} & \multicolumn{7}{c|}{$\bbeb{\mclr}$ of Kernel Matrix $K_{edge}$}  \\
            \cline{2-8}
            & $n=225$ & $n=484$ & $n=961$ & $n=1936$  & $n=3969$& $n=8100$& $n=16129$\\
            \hline
            $\mclk_1$ &80.43 &114.19 &154.62 &210.49 &292.42 &404.35 &545.65\\ 
            $\mclk_2$ &46.99 &62.52 &81.04 &106.71 &143.78 &195.15 &258.36\\
            $\mclk_3$ &69.63 &94.42 &121.97 &159.32 &208.72 &273.59 &352.50\\
            $\mclk_4$ &80.94 &114.56 &155.03 &211.29 &292.82 &405.00 &562.63\\
            $\mclk_5$ &65.74 &89.04 &115.61 &151.20 &198.90 &261.39 &337.23 \\
            $\mclk_6$ &69.84 &95.16 &124.15 &161.81 &215.79 &284.33 &368.51\\
            $\mclk_7$ &67.60 &91.14 &118.67 &155.42 &203.42 &266.21 &343.58\\
            \hline
        \end{tabular}}
        \subcaption{Mean of Numerical Ranks of Data Samples}
        \label{tab: expected ranks edge sharing domain in 2D}
    \end{minipage}%
    \hfill 
    \begin{minipage}{0.49\textwidth}
        \centering
        \resizebox{\textwidth}{!}{ 
        \begin{tabular}{|c| c c c c c c c|}
            \hline
            \multirow{2}{*}{ker fun} & \multicolumn{7}{c|}{$\var{\mclr}$ of Kernel Matrix $K_{edge}$}  \\
            \cline{2-8}
                & $n=225$ & $n=484$ & $n=961$ & $n=1936$  & $n=3969$& $n=8100$& $n=16129$\\
            \hline
            $\mclk_1$ &111.77 &259.58 &542.42 &1029.01 &2211.46 &4461.50 &9999.78\\ 
            $\mclk_2$ &63.84 &127.42 &256.70 &497.60 &1038.78 &2097.10 &4340.36\\
            $\mclk_3$ &90.48 &174.31 &327.43 &624.99 &1182.06 &2174.75 &4698.81\\
            $\mclk_4$ &110.98 &257.50 &516.83 &1046.45 &2228.99 &4421.87 &9131.77\\
            $\mclk_5$ &84.96 &174.31 &327.43 &624.99 &1182.06 &2174.75 &4178.44 \\
            $\mclk_6$ &89.11 &167.77 &381.23 &719.93 &1284.90 &2358.76 &5048.11\\
            $\mclk_7$ &89.79 &186.05 &349.42 &636.91 &1198.51 &2239.94 &4475.94\\
            \hline
        \end{tabular}}
        \subcaption{Variance of Numerical Ranks of Data Samples}
        \label{tab: variance ranks edge sharing domain in 2D}
    \end{minipage}
    \caption{ Random rank statistics for edge-sharing domains in 2D.}
\end{table}


\subsubsection*{Results for Three-Dimensional Cases:}

\begin{table}[H]
    \centering
    \begin{minipage}{0.49\textwidth}
        \centering
        \resizebox{\textwidth}{!}{ 
        \begin{tabular}{|c| c c c c c c c|}
            \hline
            \multirow{2}{*}{ker fun} & \multicolumn{7}{c|}{$\bbeb{\mclr}$ of Kernel Matrix $K_{far}$}  \\
            \cline{2-8}
                & $n=216$ & $n=512$ & $n=1331$ & $n=2197$  & $n=4096$& $n=8000$& $n=15625$\\
            \hline
            $\mclk_1$ &42.71 &47.47 &49.67 &50.50 &51.33 &52.48 &53.77\\ 
            $\mclk_2$ &47.75 &54.60 &60.88 &63.39 &65.09 &66.07 & 66.70\\
            $\mclk_3$ &60.05 &68.84 &74.92 &77.19 &78.92 &80.52 &81.90\\
            $\mclk_4$ &45.31 &49.22 &52.53 &54.23 &56.66 &59.21 &61.16\\
            $\mclk_5$ &45.63 &50.91 &56.28 &58.83 &61.38 &63.50 &65.05 \\
            $\mclk_6$ &55.71 &63.89 &69.89 &71.92 &73.62 &74.9 &75.95\\
            $\mclk_7$ &44.94 &49.55 &53.50 &55.40 &57.05 &59.02 &60.83\\
            \hline
        \end{tabular}}
        \subcaption{Mean of Numerical Ranks of Data Samples}
        \label{tab: expected ranks far-field domain in 3D}
    \end{minipage}%
    \hfill 
    \begin{minipage}{0.49\textwidth}
        \centering
        \resizebox{\textwidth}{!}{ 
        \begin{tabular}{|c| c c c c c c c|}
            \hline
            \multirow{2}{*}{ker fun} & \multicolumn{7}{c|}{$\var{\mclr}$ of Kernel Matrix $K_{far}$}  \\
            \cline{2-8}
                & $n=216$ & $n=512$ & $n=1331$ & $n=2197$  & $n=4096$& $n=8000$& $n=15625$\\
            \hline
            $\mclk_1$ &20.88 &9.33 &3.93 &4.71 &6.51 &8.00 &9.05\\ 
            $\mclk_2$ &19.30 &31.62 &23.35 &12.88 &5.05 &1.59 &1.02\\
            $\mclk_3$ &33.20 &30.75 &16.46 &11.16 &8.23 &6.07 &3.86\\
            $\mclk_4$ &18.31 &6.25 &16.11 &20.34 &21.90 &17.97 &8.88\\
            $\mclk_5$ &13.78 &19.01 &23.04 &20.71 &15.37 &8.72 &3.26 \\
            $\mclk_6$ &24.97 &26.95 &14.15 &9.07 &5.29 &4.10 &3.37\\
            $\mclk_7$ &12.89 &11.82 &12.32 &12.99 &12.69 &10.56 &8.31\\
            \hline
        \end{tabular}}
        \subcaption{Variance of Numerical Ranks of Data Samples}
        \label{tab: variance ranks far-field domain in 3D}
    \end{minipage}
    \caption{ Random rank statistics for far-field domains in 3D.}
\end{table}

\begin{table}[H]
    \centering
    \begin{minipage}{0.49\textwidth}
        \centering
        \resizebox{\textwidth}{!}{ 
        \begin{tabular}{|c| c c c c c c c|}
            \hline
            \multirow{2}{*}{ker fun} & \multicolumn{7}{c|}{$\bbeb{\mclr}$ of Kernel Matrix $K_{ver}$}  \\
            \cline{2-8}
                & $n=216$ & $n=512$ & $n=1331$ & $n=2197$  & $n=4096$& $n=8000$& $n=15625$\\
            \hline
            $\mclk_1$ &80.38 &102.51 &127.34 &140.71 &155.51 &171.64 &188.24\\ 
            $\mclk_2$ &95.44 &127.39 &162.17 &179.04 &198.18 &219.72 &239.93\\
            $\mclk_3$ &92.43 &118.24 &144.71 &157.15 &171.73 &186.32 &200.58\\
            $\mclk_4$ &81.73 &105.05 &129.85 &142.17 &158.36 &174.77 &188.71\\
            $\mclk_5$ &84.13 &109.08 &135.07 &147.42 &160.94 &175.72 &189.85 \\
            $\mclk_6$ &93.00 &122.25 &151.82 &166.75 &183.09 &200.73 &217.53\\
            $\mclk_7$ &83.78 &106.68 &131.00 &142.40 &155.01 &169.00 &181.81\\
            \hline
        \end{tabular}}
        \subcaption{Mean of Numerical Ranks of Data Samples}
        \label{tab: expected ranks vertex sharing domain in 3D}
    \end{minipage}%
    \hfill 
    \begin{minipage}{0.49\textwidth}
        \centering
        \resizebox{\textwidth}{!}{ 
        \begin{tabular}{|c| c c c c c c c|}
            \hline
            \multirow{2}{*}{ker fun} & \multicolumn{7}{c|}{$\var{\mclr}$ of Kernel Matrix $K_{ver}$}  \\
            \cline{2-8}
                & $n=216$ & $n=512$ & $n=1331$ & $n=2197$  & $n=4096$& $n=8000$& $n=15625$\\
            \hline
            $\mclk_1$ &114.36 &193.29 &276.27 &322.43 &368.38 &391.90 &410.98\\ 
            $\mclk_2$ &172.11 &348.35 &484.78 &528.45 &654.80 &642.32 &688.23\\
            $\mclk_3$ &129.73 &223.50 &277.14 &322.66 &339.84 &337.85 &287.35\\
            $\mclk_4$ &119.24 &181.67 &277.96 &313.00 &345.77 &373.36 &515.64\\
            $\mclk_5$ &130.20 &214.16 &279.47 &298.08 &321.24 &335.96 &327.66\\
            $\mclk_6$ &130.19 &229.05 &327.38 &376.41 &380.03 &397.55 &377.13\\
            $\mclk_7$ &128.02 &200.02 &272.62 &303.72 &312.34 &323.02 &266.42\\
            \hline
        \end{tabular}}
        \subcaption{Variance of Numerical Ranks of Data Samples}
        \label{tab: variance ranks ver sharing domain in 3D}
    \end{minipage}
    \caption{ Random rank statistics for vertex-sharing domains in 3D.}
\end{table}

\begin{table}[H]
    \centering
    \begin{minipage}{0.49\textwidth}
        \centering
        \resizebox{\textwidth}{!}{ 
        \begin{tabular}{|c| c c c c c c c|}
            \hline
            \multirow{2}{*}{ker fun} & \multicolumn{7}{c|}{$\bbeb{\mclr}$ of Kernel Matrix $K_{edge}$}  \\
            \cline{2-8}
                & $n=216$ & $n=512$ & $n=1331$ & $n=2197$  & $n=4096$& $n=8000$& $n=15625$\\
            \hline
            $\mclk_1$ &97.14 &135.99 &186.18 &215.26 &258.26 &313.25 &377.63\\ 
            $\mclk_2$ &114.98 &166.89 &234.70 &277.50 &332.52 &399.49 &487.64\\
            $\mclk_3$ &107.24 &148.14 &198.82 &228.74 &265.06 &311.22 &363.59\\
            $\mclk_4$ &98.04 &136.73 &187.38 &218.90 &264.43 &314.35 &392.67\\
            $\mclk_5$ &99.58 &139.77 &188.84 &217.89 &256.25 &299.12 &345.54 \\
            $\mclk_6$ &108.41 &152.84 &208.12 &242.17 &283.92 &334.78 &392.67\\
            $\mclk_7$ &100.21 &137.97 &186.32 &212.34 &248.60 &290.42 &345.77\\
            \hline
        \end{tabular}}
        \subcaption{Mean of Numerical Ranks of Data Samples}
        \label{tab: expected ranks edge sharing domain in 3D}
    \end{minipage}%
    \hfill 
    \begin{minipage}{0.49\textwidth}
        \centering
        \resizebox{\textwidth}{!}{ 
        \begin{tabular}{|c| c c c c c c c|}
            \hline
            \multirow{2}{*}{ker fun} & \multicolumn{7}{c|}{$\var{\mclr}$ of Kernel Matrix $K_{edge}$}  \\
            \cline{2-8}
                & $n=216$ & $n=512$ & $n=1331$ & $n=2197$  & $n=4096$& $n=8000$& $n=15625$\\
            \hline
            $\mclk_1$ &191.54 &459.26 &968.90 &1341.32 &2068.00 &3127.74 &4928.74\\ 
            $\mclk_2$ &266.14 &719.88 &1618.99 &2219.04 &3171.97 &4821.65 &6527.24\\
            $\mclk_3$ &210.26 &466.79 &935.47 &1348.08 &1723.51 &2516.34 &3166.02\\
            $\mclk_4$ &202.77 &458.58 &968.32 &1337.33 &2118.39 &3212.67 &3844.92\\
            $\mclk_5$ &200.48 &480.89 &929.38 &1231.47 &1778.18 &2278.84 &3075.45 \\
            $\mclk_6$ &200.48 &476.91 &1017.19 &1393.06 &1944.73 &2931.39 &3844.92\\
            $\mclk_7$ &210.49 &438.54 &919.84 &1235.32 &1655.63 &2356.37 &3063.37\\
            \hline
        \end{tabular}}
        \subcaption{Variance of Numerical Ranks of Data Samples}
        \label{tab: variance ranks edge sharing domain in 3D}
    \end{minipage}
    \caption{ Random rank statistics for edge-sharing domains in 3D.}
\end{table}


\begin{table}[H]
    \centering
    \begin{minipage}{0.49\textwidth}
        \centering
        \resizebox{\textwidth}{!}{ 
        \begin{tabular}{|c| c c c c c c c|}
            \hline
            \multirow{2}{*}{ker fun} & \multicolumn{7}{c|}{$\bbeb{\mclr}$ of Kernel Matrix $K_{face}$}  \\
            \cline{2-8}
                & $n=216$ & $n=512$ & $n=1331$ & $n=2197$  & $n=4096$& $n=8000$& $n=15625$\\
            \hline
            $\mclk_1$ &126.05 &204.18 &342.51 &451.14 &637.05 &925.74 &1358.29\\ 
            $\mclk_2$ &145.12 &244.05 &416.10 &544.89 &766.05 &1103.68 &1620.39\\
            $\mclk_3$ &132.68 &213.20 &345.32 &444.86 &607.73 &855.64 &1195.01\\
            $\mclk_4$ &126.29 &207.91 &344.50 &449.07 &639.83 &937.05 &1329.50\\
            $\mclk_5$ &126.33 &203.83 &331.93 &429.36 &581.88 &819.84 &1177.91\\
            $\mclk_6$ &133.89 &218.05 &355.14 &457.75 &627.96 &888.22 &1248.62\\
            $\mclk_7$ &128.53 &205.55 &334.70 &429.56 &588.33 &825.41 &1183.07\\
            \hline
        \end{tabular}}
        \subcaption{Mean of Numerical Ranks of Data Samples}
        \label{tab: expected ranks face sharing domain in 3D}
    \end{minipage}%
    \hfill 
    \begin{minipage}{0.49\textwidth}
        \centering
        \resizebox{\textwidth}{!}{ 
        \begin{tabular}{|c| c c c c c c c|}
            \hline
            \multirow{2}{*}{ker fun} & \multicolumn{7}{c|}{$\var{\mclr}$ of Kernel Matrix $K_{face}$}  \\
            \cline{2-8}
                & $n=216$ & $n=512$ & $n=1331$ & $n=2197$  & $n=4096$& $n=8000$& $n=15625$\\
            \hline
            $\mclk_1$ &335.66 &1399.68 &5227.05 &9623.13 &23369.13 &55586.36 &117362.16\\ 
            $\mclk_2$ &335.89 &1490.79 &5721.25 &12238.55 &25217.08 &58863.42 &137613.76\\
            $\mclk_3$ &326.53 &1293.61 &4152.90 &8770.68 &18930.69 &44508.93 &82545.65\\
            $\mclk_4$ &313.89 &1288.67 &4894.75 &10180.30 &22019.31 &54670.16 &117559.78\\
            $\mclk_5$ &318.03 &1216.24 &4121.56 &8162.81 &17599.36 &36474.29 &89411.12\\
            $\mclk_6$ &304.99 &1233.65 &4788.07 &8644.32 &19568.32 &44188.66 &95840.63\\
            $\mclk_7$ &337.33 &1264.44 &4505.36 &8595.78 &17812.45 &40953.21 &88895.33\\
            \hline
        \end{tabular}}
        \subcaption{Variance of Numerical Ranks of Data Samples}
        \label{tab: variance ranks face sharing domain in 3D}
    \end{minipage}
    \caption{ Random rank statistics for face-sharing domains in 3D.}
\end{table}

\section{Conclusion}
\label{sec: Conclusion}
In this article, we explored the behavior of the rank of the kernel matrices arising from the interactions of neighboring source and target domains, where particles are arbitrarily distributed, moving beyond the typical assumption of uniform grid settings, and the arbitrary distribution of particles was modeled to arise from an underlying random distribution. We have established theoretical results on the growth of the expectation and variance of the random rank $\mclr$ for all possible neighboring interactions in $d$-dimensions. The numerical results in one-, two-, and three-dimensions confirmed our theoretical predictions as stated in \autoref{sec: Main Results}. Our analysis guarantees that despite the inherent arbitrariness in particle distributions, the rank structure of matrices due to all possible interactions between neighboring source and target domains remains consistent, which offers a strong guarantee for the efficiency of hierarchical matrix algorithms in real-world applications. Moreover, the approach of choosing randomly distributed particles allows us to examine algorithmic behavior in a way that mirrors the average-case analysis of algorithms, making our findings highly relevant. 

Interestingly, our observations suggest that the random variable $\mclr$, due to vertex-sharing interaction, may follow a Gaussian distribution. This raises an exciting possibility for future research, and we leave the formal proof as an open question.

By presenting a novel study of random rank $\mclr$ due to randomly distributed particles in the domains, this work also opens the door for more comprehensive analyses of algorithm performance in diverse practical scenarios. Future research could extend these results to more complex kernels or explore additional probabilistic frameworks for further enhancing the applicability.

\section*{Acknowledgments}
The authors would like to thank the High-Performance Computing Environment (HPCE) at IIT Madras for providing the computational resources essential to this work. Additionally, the first author acknowledges the financial support from UGC through the Junior Research Fellowship (JRF) awarded for the doctoral research and also would like to thank the IIT Madras Library for providing access to Grammarly, which greatly assisted in improving the grammar of the article.

{\footnotesize
\bibliographystyle{siam}
\bibliography{ref}}

\appendix

\section{Generalization to Arbitrary Probability Distributions} \label{app: Generalization to Arbitrary Probability Distributions}

We have already discussed in \autoref{subsec: Choice of Probability Distribution} about choosing the uniform probability distribution in detail. Now, in this section, we give an idea of how to generalize the probability distribution rather than sticking only to the uniform probability distribution, which allows for more flexible probabilistic modeling. For that, let us consider $ x_1, x_2, \dotsc, x_n $ be $n$ i.i.d. random particles in $ V= [a,b]^d $. 

\subsection{Identical Marginals with Coordinate-wise Independence:}
Let us assume that each particle $ x_i = \bkt{x_i^{(1)}, x_i^{(2)}, \dotsc, x_i^{(d)}} $ in $V$ is drawn from some \textit{arbitrary product distribution with i.i.d. coordinates}. That is, there exists a univariate probability distribution with CDF $ \psi $, supported on $[a,b]$, such that $x_i^{(j)} \sim \psi, \text{ for all } j = 1, \dotsc, d.$ Let \( V' = \prod_{j=1}^d[c_j,d_j] \subseteq V \) be a sub-hyper-cube, and the random variable $ N(V') $ is defined as the number of particles that fall within $V'$. Then the probability that a single particle lies within $V'$ is
\[
\mathtt{q}_{_{V'}} = \prod_{j=1}^d \bkt{ \frac{\psi(d_j) - \psi(c_j)}{\psi(b) - \psi(a)} },
\]
and the number of particles that fall within $V'$ follows the Binomial distribution as given by \autoref{equ: binomial distribution pmf}.

\subsection{Non-identical Marginals but Coordinate-wise Independence:}
More generally, let us consider that each particle $ x_i = \bkt{x_i^{(1)}, x_i^{(2)}, \dotsc, x_i^{(d)}} $ in $V$ is drawn from some \textit{arbitrary product distribution}, but each coordinate \( x_i^{(j)} \) follows \textit{different distribution} with CDF \( \psi_j \), supported on \( [a,b] \), i.e., $ x_i^{(j)} \sim \psi_j, \text{ for each } j = 1, \dotsc, d $, and hence, the probability that a single particle lies within $V'$ is
\[
\mathtt{q}_{_{V'}} = \prod_{j=1}^{d} \bkt{\frac{\psi_j(d_j) - \psi_j(c_j)}{\psi_j(b) - \psi_j(a)}}.
\]
Accordingly, the number of particles that fall within $V'$ follows the Binomial distribution as given by \autoref{equ: binomial distribution pmf}.

\begin{myremark}
    The assumption that the particles in $V$ are coordinate-wise independent is essential for both cases. Without it, the joint distribution of $ x_i $ may not admit a product structure, and the distribution of \( N(V') \) may not be Binomial.
\end{myremark}


\section{Calculation of \texorpdfstring{ $\bbeb{N(V')N(V'')}$}{E[N(V')N(V'')]} and \texorpdfstring{$ \cov{N(V'),N(V'')} $}{Cov(N(V')N(V''))}} \label{app: Calculation of E[N(V')N(V'')] and Cov(N(V'),N(V''))}
Here, $V'$ and $V''$ are two non-intersecting subdomains of $V$, and the random variables $N(V')$ and $N(V'')$ count the number of particles in $V'$ and $V''$, respectively. Also, $\mathtt{q}_1$ is the probability of having a particle in $V'$, and that of $V''$ is $\mathtt{q}_2$. Now, we define following two indicator random variables \[ \mathbb{I}'_i=\begin{cases} 1&\text{if }x_i\in V'\\0 &\text{otherwise} \end{cases}\quad \text{and}\quad \mathbb{I}''_i=\begin{cases} 1&\text{if }x_i\in V''\\0 &\text{otherwise} \end{cases}. \]  As $ \displaystyle N(V')N(V'')= \sum_{i=1}^n\sum_{j=1}^n \mathbb{I}'_i\mathbb{I}''_j $, we have $ \displaystyle \bbeb{N(V')N(V'')}= \sum_{i=1}^n\sum_{j=1}^n \bbeb{\mathbb{I}'_i\mathbb{I}''_j} $. Now, we have the following cases \begin{itemize}
    \item If $i=j$, then $\bbpb{x_i\in V'\text{ and }x_i\in V''}=0$, as the domains are non-intersecting. This gives $\bbeb{\mathbb{I}'_i\mathbb{I}''_i}=0$.
    \item If $i\neq j$, then we have $\bbeb{\mathbb{I}'_i\mathbb{I}''_j} = \bbeb{\mathbb{I}'_i}\bbeb{\mathbb{I}''_j} = \mathtt{q}_1\mathtt{q}_2$, as the particles $x_i$'s are i.i.d. in $V$.
\end{itemize} Thus we have $ \displaystyle \bbeb{N(V')N(V'')} = \sum_{i=1}^n\sum_{j=1}^n \bbeb{\mathbb{I}'_i}\bbeb{\mathbb{I}''_j} = \sum_{i=1}^n\sum\limits_{\substack{j=1 \\ j \neq i}}^n \mathtt{q}_1\mathtt{q}_2 = n(n-1)\mathtt{q}_1\mathtt{q}_2 $. Now, $ \cov{N(V'),N(V'')} = -n\mathtt{q}_1\mathtt{q}_2 $ follows directly from the previous expectation and the fact that $N(V')$ and $N(V'')$ follow binomial distribution individually.

\section{Relationship between numerical \texorpdfstring{$\varepsilon$}{epsilon}-rank and max-rank} \label{app: Relationship between numerical rank and max-rank}
This is a well-known fact that $\magn{\cdot}_{\infty^*}$ and $\magn{\cdot}_2$ are equivalent norms and the corresponding equivalency relation is \begin{equation} \label{equ: norm equivalency}
    \magn{A}_{\infty^*} \leq \magn{A}_2 \leq \sqrt{mn}\magn{A}_\infty, \quad \text{where } A\in \bbc^{m\times n}.
\end{equation} 
Now, to establish the relationship between numerical $\varepsilon$-rank and max-rank, it is sufficient to show how $\displaystyle \frac{\magn{A-\tilde{A}}_{\infty^*}}{\magn{A}_{\infty^*}}$ and $\displaystyle \frac{\magn{A-\tilde{A}}_2}{\magn{A}_2}$ are related, where $\Tilde{A}$ is the approximation of $A$. From \autoref{equ: norm equivalency}, we can easily get the following \[ \frac{1}{\sqrt{mn}} \frac{\magn{A-\tilde{A}}_2}{\magn{A}_2} \leq \frac{\magn{A-\tilde{A}}_{\infty^*}}{\magn{A}_{\infty^*}} \leq \sqrt{mn} \frac{\magn{A-\tilde{A}}_2}{\magn{A}_2} \] Now from the above expression we have the followings: \begin{itemize}
    \item For any given $\varepsilon>0$, there exist $\varepsilon'>0$ such that \[ \frac{\magn{A-\tilde{A}}_2}{\magn{A}_2} < \varepsilon \implies \frac{\magn{A-\tilde{A}}_{\infty^*}}{\magn{A}_{\infty^*}} <\varepsilon', \quad\text{where } \varepsilon' = \sqrt{mn}\varepsilon. \]
    \item For any given $\varepsilon'>0$, there exist $\varepsilon>0$ such that \[ \frac{\magn{A-\tilde{A}}_{\infty^*}}{\magn{A}_{\infty^*}} < \varepsilon' \implies \frac{\magn{A-\tilde{A}}_2}{\magn{A}_2} <\varepsilon, \quad\text{where } \varepsilon = \sqrt{mn}\varepsilon'. \]
\end{itemize}

\section{Error Due to Normal Approximation} \label{app: Error Due to Normal Approximation}
 In this section, we will focus on the error bounds for normal approximations in both the one-dimensional and $d$-dimensional cases. The Berry-Esseen theorem \cite{Feller1971AnIntro}, and its extensions to higher dimensions \cite{Bentkuus2004ALyap} provide a framework to quantify the error in the normal approximation to the sum of independent random variables. These error bounds provide an understanding of the rate of convergence to the normal distribution, which is useful in practical applications.
\subsection{One-dimensional case:}
Let $X_1,X_2,\dotsc, X_n$ be i.i.d random variables with $\bbeb{X_k}=0$, $\bbeb{X_k^2}=\sigma^2>0$ and $\bbeb{X_k^3}=\rho<\infty$. Then the cumulative distribution function $F_n$ of the normalized sum $S_n = \dfrac{X_1+X_2+\dotsb+X_n}{\sigma\sqrt{n}}$ converges to the standard normal distribution $\Phi$ with \[ \abslt{F_n(x)-\Phi(x)}\leq \frac{3\rho}{\sigma^3\sqrt{n}}\quad\quad \forall\; x,n. \]
The proof can be found in \cite{Feller1971AnIntro}. Now, as the binomial distribution can be interpreted as the sum of independent Bernoulli's and hence from the above expression, we can get an error bound in the normal approximation to the binomial distribution.
\subsubsection*{Application to Binomial Distribution}
Suppose $X_1,X_2,\dotsc,X_n$ are i.i.d. Bernoulli random variables with success probability $p$, and $S_n=\sum\limits_{i=1}^n X_i $ is a Binomial random variable, i.e., $S_n \sim \text{Bin}(n,p)$ \cite{Lebanon2013Probability}. Now let \( X \sim \text{Bernoulli}(p) \), where \( X \) takes the value 1 with probability \( p \) and 0 with probability \( 1-p \). Then, standard deviation $\sigma = \sqrt{\var{X}} = \sqrt{p(1-p)}$ and the third absolute central moment $\rho = \bbeb{|X - \mu|^3} = p(1-p)(1 - 2p + 2p^2)$. Now, using the Berry-Esseen theorem, the error in approximating the binomial distribution with a normal distribution is:
\begin{equation}
    \sup_x \left| P\left( \frac{S_n - np}{\sqrt{np(1-p)}} \leq x \right) - \Phi(x) \right| \leq \frac{ C (1-2p+2p^2)}{\sqrt{n p (1-p)}}.
\end{equation}

\subsection{\texorpdfstring{$d$}{d}-dimensional case:}
Let $X_1,X_2,\dotsc,X_n$ be independent and identically distributed random vectors in $\bbr^d$ such that for each $k = 1:n$, $\bbeb{X_k}=0$ and $X_k$ has identity covariance matrix. Now let $S_n=X_1+X_2+\dotsb+X_n$ and $Z$ be the Gaussian random variable with the same mean and variance that of $S_n$. Now for any convex subset $\mathcal{C}$ in $\bbr^d$, we have a Lyapunov-type bound as \[ \sup_{A\in \mathcal{C}}\abslt{\bbpb{S_n\in A} - \bbpb{Z\in A} } \leq \frac{cd^{1/4}\bbeb{\magn{X_1}_2^3}}{\sqrt{n}}, \text{ where $c$ is a constant. } \] 
The proof can be found here \cite{Bentkuus2004ALyap}. Similar to the univariate case, in the multivariate case, the multinomial distribution can be interpreted as the sum of independent multivariate Bernoulli's (or simply Multinoulli) \cite{Lebanon2013Probability} and hence from the above statement, we can get the error due to the multivariate normal approximation to the multinomial distribution (i.e., for $d=2$) 

\subsubsection*{Application to Multinomial}
Suppose $X_1, X_2,\dotsc, X_n$ are i.i.d. Multivariate Bernoulli random variables with success probability $p = (p_1,p_2,\dotsc,p_d)$, and $S_n=\sum\limits_{i=1}^n X_i $ is a Multinomial random variable, i.e., $S_n \sim \text{Mult}(n,p)$ \cite{Lebanon2013Probability}. Now let \( X \sim \text{Multinoulli}(p) \), i.e., $X$ is a set of $d\times 1$ vectors having one entry equal to 1 and all other entries equal to 0 and $p_1,p_2,\dotsc,p_d$ such that $0<p_i<1$, for each $i$ with the condition that $\displaystyle\sum_{i=1}^dp_i=1$. Then, the joint probability mass function of $X$ is given by \[ p_X(x_1,x_2,\dotsc,x_d)= \begin{cases}
    \prod_{k=1}^d p_k^{x_j} &\text{if } (x_1,x_2,\dotsc,x_d)^t\in X \\ 0 &\text{otherwise.}
\end{cases} \] i.e., $X$ takes value $e_i$ with probability $p_i$. Now, we have $\bbeb{X} = p$, $\Sigma = \bkt{\Sigma_{ij}}_{d\times d}$  the covariance matrix of $X$, where $ \Sigma_{ij} = \begin{cases} p_i(1-p_i) &\text{if } i=j \\ -p_ip_j &\text{if } i\neq j \end{cases}$ and as for all $i$, $\magn{X-p}_2^2 = \magn{e_i-p}_2^2 = 1 - 2p_i + \sum_{k=1}^dp_k^2 $, then \[ \bbeb{ \magn{X-p}_2^3 } = \sum_{i=1}^d p_i \bkt{1 - 2p_i + \sum_{k=1}^dp_k^2}^{{3}/{2}} < \frac{M}{cd^{{1}/{4}}} \quad \text{ (a constant) (say)}. \]
Now, to get the error-bound while approximating $S_n$ with $Z\sim\mathcal{N}\bkt{np,n\Sigma}$, we transform $S_n$ to $W_n$ such that $W_n = \Sigma^{-1/2}\bkt{S_n -np} $ and similarly $Z$ to $Z'$, where $Z'\sim \mathcal{N}\bkt{0,I_d}$. Then we have  \[ \sup_{A\in \mathcal{C}}\abslt{\bbpb{S_n\in A} - \bbpb{Z\in A} } = \sup_{A'\in \mathcal{C}}\abslt{\bbpb{W_n\in A'} - \bbpb{Z'\in A'} } \leq \frac{M}{\sqrt{n}}, \text{ where } A'= \Sigma^{-1/2}\bkt{A - np}. \]

\section{Uniform Distribution as Randomly Uniform Perturbed Grid}\label{app: Uniform Distribution as Randomly Uniform Perturbed Grid}
In computational applications, we often use uniform grids to discretize a continuous domain. However, in practice, random variations often lead to deviations from perfect uniformity. The perspective of taking Uniform Probability distribution allows us to interpret it as a `perturbed' version of a perfect uniform grid. For the sake of simplicity, we will consider the domain as $\mathcal{D}=[0,1]$.

Suppose that $X_1, X_2, \dotsc, X_n$ be a random sample of size $n$ from a uniform distribution $\mathcal{U}_{\mathcal{D}}$ with pdf $f$ and CDF $F$ and $X_{(1)}, X_{(2)},\dotsc, X_{(n)}$ be the corresponding order statistics, then the pdf of the $k^{th}$ order statistic is given as \[ f_{X_{(k)}}(x) = \dfrac{n!}{(k-1)! (n-k)!} f(x)[F(x)]^{k-1}[1-F(x)]^{n-k} \] 
Now, as the $k^{th}$ order statistic in a sample of size $n$ from $\mathcal{U}_{\mathcal{D}}$ has a beta distribution \cite{gentle2009computational} with parameters $k$ and $n-k+1$, the expected value and variance of $X_{(k)}$ are \[ \bbeb{X_{(k)}}=\frac{k}{n+1}\quad \text{and}\quad \var{X_{(k)}}=\frac{k\bkt{n-k+1}}{\bkt{n+1}^2\bkt{n+2}} \]
Thus, for large $n$, the expected value of $X_{(k)}$ is very close to the ideal $k^{th}$ grid point. In addition to the expectation, for large $n$, this variance becomes small, indicating that each $X_{(k)}$ tends to stay close to the ideal grid point. Moreover, if $G_k = X_{(k+1)} - X_{(k)} $ be the gap between two consecutive shorted random points, then $\bbeb{G_k}=\dfrac{1}{n+1}$, which suggest that as $n$ grows, sorted uniform points tend to fill up the domain $\mathcal{D}$ in a manner that resembles a uniform grid with smaller fluctuations.

\section{Detailed Proofs of \autoref{lma: sum of Znk1s is constant} and \autoref{lma: sum of ZniZnjs is constant}}
\label{app: proof of lemma: sum of Zni1s is constant and lemma: sum of ZniZnjs is constant}
We provide a comprehensive proof of \autoref{lma: sum of Znk1s is constant} and \autoref{lma: sum of ZniZnjs is constant}, which are briefly justified in \autoref{sec: Fundamental Framework and Problem Setup} of the main body of the article. We now restate the lemma below for clarity before the detailed proof. \vskip 0.15cm
\noindent\textbf{\autoref{lma: sum of Znk1s is constant} (Restated): } \textit{For any large value of $n$, if $\kappa = \lfloor \log_{2^d}n\rfloor $ then for some fixed $l$ the sum of variances $\displaystyle \sum_{k=1}^\kappa \var{\mclz_n^{k,l}} $ is bounded by $C$ such that $C\in \mclo{\log\log\log n}$.}
\begin{proof}  
    As \autoref{equ: variance Znkl with probability} is independent of $l$, without loss of generality, we take $l=1$.
    Now, we have \[\var{\mclz_n^{k,1}}\leq \frac{p(p+1)(2p+1)}{6}\bbp\bkt{N_{k,1}=p}\text{ and }\bbp\bkt{N_{k,1}=p} \leq \frac{1}{\sqrt{2\pi}}e^{-{\bkt{a_{p,k}^{(n)}}^{2}}/{2}} + \frac{ \Xi }{\sqrt{n}} ,\] where $a_{p,k}^{(n)}=\frac{p-0.5-\mu_k}{\sigma_k}$ with $\mu_k=\frac{n}{2^{dk}}$ and $\sigma_k=\sqrt{\frac{n}{2^{dk}}\bkt{1-\frac{1}{2^{dk}}}} $. Now, for any small $\varepsilon >0$, there exists $\Tilde{k}$ such that $\bbp\bkt{N_{k,1}=p} < \varepsilon$ for all $k\leq \Tilde{k}$ and $\Tilde{k}$ can be determined as follows
    \begin{align*}
        & \frac{1}{\sqrt{2\pi}}e^{-{\bkt{a_{p,k}^{(n)}}^{2}}/{2}} + \frac{p(p+1)}{2}\frac{ \Xi }{\sqrt{n}} <\varepsilon\\
        \implies&\frac{1}{\sqrt{2\pi}}e^{-\bkt{a_{p,k}^{(n)}}^2/2} < \varepsilon \\
        \implies &a_{p,k}^{(n)} > \sqrt{2\log_e{\frac{1}{\varepsilon\sqrt{2\pi}}}} =: M(\varepsilon) \text{ (say)}\\
        \implies & \frac{p-0.5-\frac{n}{2^{dk}}}{\sqrt{\frac{n}{2^{dk}}\bkt{1-\frac{1}{2^{dk}}}}} > M\\
        \implies & \bkt{M^2+n}nx^2 -\bkt{2\bkt{p-0.5}+M^2}nx +\bkt{p-0.5}^2 >0, \hskip 0.5cm \bkbt{\text{where } x = \frac{1}{2^{dk}}.}
    \end{align*}
    Now for the sake of simplicity, we take $\bkt{M^2+n}x -\bkt{2\bkt{p-0.5}+M^2} >0$ which implies 
    \begin{align*}
          k & < \log_{2^d}\bkt{1 + \frac{n+1-2p}{M^2 +2p-1}}
    \end{align*}
    Now we can take $\Tilde{k}= \left\lfloor \log_{2^d}\bkt{1 + \frac{n+1-2p}{M^2 +2p-1}} \right\rfloor$, then for all $k=1:\Tilde{k}$ (in MATLAB notation) we have \[\bbp\bkt{N_{k,1}=p}<\varepsilon.\]    
    There are only constantly many terms are there from $\Tilde{k}$ to $\kappa$, and the constant depends on $\varepsilon$ only which can be shown as 
    \begin{align*}
        \kappa - \Tilde{k} 
            & \leq 1+ \left\lfloor \log_{2^d}n - \log_{2^d}\bkt{ \frac{n+1-2p}{M^2 +2p-1}} \right\rfloor \\
            &\leq 1 + \log_{2^d}\bkt{{M^2 +2p-1} } + \log_{2^d}\bkt{\frac{n}{n+1-2p}} \\
            &< \Omega + \log_{2^d}\bkt{{M^2 +2p-1} }, \text{ where $\Omega$ is independent of $n$.}
    \end{align*}
    Now, we choose $\varepsilon = \frac{6c_1}{\kappa p(p+1)(2p+1)} $, where $c_1$ is a constant such that 
    \begin{align*}
        \sum_{k=1}^{\Tilde{k}} \var{ \mclz_n^{k,1} } \leq \sum_{k=1}^{\Tilde{k}} \frac{p(p+1)(2p+1)}{6} \bbp\bkt{N_{k,1}=p} < \frac{p(p+1)(2p+1) \kappa \varepsilon }{6} = c_1
    \end{align*}
    Now, 
    \begin{align*}
        \sum_{k=1}^\kappa \var{ \mclz_n^{k,1} } &=\sum_{k=1}^{\Tilde{k}} \var{ \mclz_n^{k,1} } + \sum_{k={\Tilde{k}}+1}^\kappa \var{ \mclz_n^{k,1} } \\
                    &<  c_1 + \sum_{k={\Tilde{k}}+1}^\kappa \var{ \mclz_n^{k,1} } 
    \end{align*}
    Now, from \autoref{equ: variance Znkl with probability}, it follows that $\var{ \mclz_n^{k,1} }\leq p(p+1)(2p+1)/6$ and  $\bkt{\kappa - \Tilde{k}}$ is of $\mclo{\log\log\log n}$ proves that there exists some $C > 0$ such that \[\displaystyle\sum_{k=1}^\kappa \var{ \mclz_n^{k,1} } < C,\quad \text{ where $C\in \mclo{\log\log\log n}$} .\]
\end{proof}

\noindent\textbf{\autoref{lma: sum of ZniZnjs is constant} (Restated): } \textit{For any large value of $n$, if $\kappa = \lfloor \log_{2^d}n\rfloor $ then for some fixed $l$ the sum of covariances \[\sum_{k_1=1}^\kappa\sum_{k_2=k_1+1}^\kappa\cov{\mclz_n^{k_1,l},\mclz_n^{k_2,l}}\leq C, \quad \text{ where $C \in \mclo{\bkt{\log\log\log n}^2}$.}\]}

\begin{proof}
    Without loss of generality, let us take $l=1$. Now for any fixed $i,j$, an upper bound of $\cov{\mclz_n^{i,1},\mclz_n^{j,1}}$ defined in \autoref{equ: cov(ZniZnj)} can be obtained as \begin{align*}
        \cov{\mclz_n^{i,1}, \mclz_n^{j,1}} \leq \sum_{r=0}^p\sum_{s=0}^p rs\bbpb{N_{i,1}=p_i,N_{j,1}=p_j}
    \end{align*} 
    for some $p_i,p_j$ with $1\leq p_i,p_j \leq p $ such that  $\bbpb{N_{i,1}=p_i,N_{j,1}=p_j} \geq \bbpb{N_{i,1}=l,N_{j,1}=m}, \forall\; 0\leq r,s\leq p $. Now, using bivariate normal distribution to approximate $\bbpb{N_{i,1}=p_i,N_{j,1}=p_j}$ we have the following \[ \bbpb{N_{i,1}=p_i,N_{j,1}=p_j}  \leq \iint\limits_{\mathcal{D}}f(x,y)dx dy + \frac{ \Xi }{\sqrt{n}}, \quad\text{ for some }  \Xi >0,  \]
    where $\mathcal{D} = \bkct{(x,y)\in\bbr^2 : a_{p_i,i}^{(n)} \leq x \leq b_{p_i,i}^{(n)} \text{ and }  a_{p_j,j}^{(n)} \leq y \leq b_{p_j,j}^{(n)}  } $ where, $a_{r,k}^{(n)} = \frac{r - 0.5 - \mu_k }{\sigma_k}$ and $b_{r,k}^{(n)} = \frac{r + 0.5 - \mu_k }{\sigma_k}$ for $r = p_i, p_j $ and $k = i,j$. 
    Now, there exist $\bkt{x_{p_i,i}^{(n)} , x_{p_j,j}^{(n)} } \in \mathcal{D} $ such that 
    \[ \iint\limits_{\mathcal{D}} f(x,y)dx dy  = f\bkt{x_{p_i,i}^{(n)} , x_{p_j,j}^{(n)}}. \] 
    Thus we have \[ \bbpb{N_{i,1}=p_i,N_{j,1}=p_j} \leq f\bkt{x_{p_i,i}^{(n)} , x_{p_j,j}^{(n)}} + \frac{ \Xi }{\sqrt{n}}. \]
    Now as $n\to \infty$, $a_{p_i,i}^{(n)} , b_{p_i,i}^{(n)}, a_{p_j,j}^{(n)}, b_{p_j,j}^{(n)} \to -\infty$ and hence $f\bkt{x_{p_i,i}^{(n)} , x_{p_j,j}^{(n)}}$ can be made sufficiently small by taking large enough $n$. Now for any chosen $\varepsilon>0$ there exist $\Tilde{k}$ such that \[ \bbpb{N_{i,1}=p_i,N_{j,1}=p_j} <\varepsilon, \quad\forall\;i,j\leq \Tilde{k}. \] 
    The largest such $\Tilde{k}$ can be found by solving the inequality  $f\bkt{x_{p_i,i}^{(n)} , x_{p_j,j}^{(n)}} + \frac{ \Xi }{\sqrt{n}} < \varepsilon $ which implies $f\bkt{x_{p_i,i}^{(n)} , x_{p_j,j}^{(n)}}<\varepsilon$. This now translates to \[ \frac{1}{2\pi\sqrt{1-\rho^2}} e^{-Q\bkt{x_{p_i,i}^{(n)} , x_{p_j,j}^{(n)}}/2\bkt{1-\rho^2}}< \varepsilon \]
    where $Q\bkt{x,y} = x^2 + y^2 - 2\rho xy $. Now further, we obtain
    \begin{align*}
        & Q\bkt{x_{p_i,i}^{(n)} , x_{p_j,j}^{(n)}} > 2\bkt{1-\rho^2}\log_e\bkt{\frac{1}{2\pi\sqrt{1-\rho^2}\varepsilon}} \\
        \implies& \bkt{x_{p_i,i}^{(n)}}^2 + \bkt{x_{p_j,j}^{(n)}}^2 -2 \rho \bkt{ x_{p_i,i}^{(n)}} \bkt{ x_{p_j,j}^{(n)}} > 2\bkt{1-\rho^2}\log_e\bkt{\frac{1}{2\pi\sqrt{1-\rho^2}\varepsilon}} 
    \end{align*}
    Now, as $\rho$ being negative, taking $x_{p_i,i}^{(n)}$ as minimum between $x_{p_i,i}^{(n)}$ and $x_{p_j,j}^{(n)}$, we have the following
    \begin{align*}
        \bkt{1-\rho}\bkt{x_{p_i,i}^{(n)}}^2 & > \bkt{1-\rho^2}\log_e\bkt{\frac{1}{2\pi\bkt{1-\rho^2}\varepsilon}} \\
        \implies \bkt{x_{p_i,i}^{(n)}}^2 & > \bkt{1+\rho}\log_e\bkt{\frac{1}{2\pi\bkt{1-\rho^2}\varepsilon}} \\
        \implies x_{p_i,i}^{(n)} & > \sqrt{ \bkt{1+\rho}\log_e\bkt{\frac{1}{2\pi\bkt{1-\rho^2}\varepsilon}} } =: M(\varepsilon);\quad\text{ (say)}\\
        \implies x_{p_i,i}^{(n)} & >  M
    \end{align*}
    Now, $x_{p_i, i}^{(n)}  >  M$ will hold also when $a_{p_i, i}^{(n)}  >  M$ is true and thus similar to the proof of previous lemma, we can prove that there exists $C>0$ such that  \[ \sum_{i=1}^\kappa\sum_{j=i+1}^\kappa\cov{\mclz_n^{i,1},\mclz_n^{j,1}} \leq C, \text{where $C\in \mclo{\bkt{\log\log\log n}^2}$}. \]
\end{proof}

\end{document}